\UseRawInputEncoding
\documentclass[11pt]{article}
\usepackage{amsmath}
\usepackage{empheq}
\usepackage{amsmath,amssymb,amsthm,graphicx,mathrsfs}
\usepackage[numbers,sort&compress]{natbib}
\usepackage{color}
\usepackage[colorlinks,bookmarksopen,bookmarksnumbered,citecolor=black, linkcolor=blue, urlcolor=black]{hyperref}
\numberwithin{equation}{section}
\newcommand{\beq}{\begin{equation}}
\newcommand{\enq}{\end{equation}}

\newtheorem{theorem}{Theorem}[section]
\newtheorem{Theorem}{Theorem}[section]
\newtheorem{Lemma}[Theorem]{Lemma}
\newtheorem{Corollary}[Theorem]{Corollary}
\newtheorem{Definition}[Theorem]{Definition}
\newtheorem{Remark}[Theorem]{Remark}

\newcommand{\benu}{\begin{enumerate}}
\newcommand{\beqa}{\begin{eqnarray}}
\newcommand{\beqan}{\begin{eqnarray*}}
\newcommand{\eay}{\end{array}}
\newcommand{\edm}{\end{displaymath}}
\newcommand{\eenu}{\end{enumerate}}
\newcommand{\eeq}{\end{equation}}
\newcommand{\eeqa}{\end{eqnarray}}
\newcommand{\eeqan}{\end{eqnarray*}}

\newcommand{\br}{\begin{Remark}}
\newcommand{\er}{\end{Remark}}

\newcommand{\bqa}{\begin{eqnarray}}
\newcommand{\eqa}{\end{eqnarray}}
\newcommand{\bqw}{\begin{eqnarray*}}
\newcommand{\eqw}{\end{eqnarray*}}

\newcommand{\non}{\nonumber}
\newcommand{\bea}{\begin{array}{cc}}
\newcommand{\ena}{\end{array}}

\allowdisplaybreaks[4]
\begin{document}
\begin{center}

{\large \bf Existence of pullback \(\mathcal{D}\)-attractors for Kirchhoff wave equations with strong damping and delay }\\

\vspace{0.20in} Bin Yang $^{1,\ast}$ $\ $ Yuming Qin $^{2}$ $\ $ Alain Miranville $^{3}$ $\ $ Xuxiao Hou $^{4}$\\
\end{center}
{\small
$^{1}$ School of Science, Inner Mongolia University of Science and Technology, Baotou 014010, Inner Mongolia, China\\
$^{2}$ School of Mathematics and Statistics, Institute for Nonlinear Science, Donghua University, Shanghai 201620, China\\
$^{3}$ Laboratoire de Math\'ematiques Appliqu\'ees du Havre (LMAH), Universit\'e Le Havre Normandie, 76063 Le Havre Cedex, France\\
$^{4}$ School of Automation and Electrical Engineering, Inner Mongolia University of Science and Technology, Baotou 014010, Inner Mongolia, China
\vspace{3mm}
}
\begin{abstract}


This paper studies the existence of pullback $\mathcal{D}$-attractors for non-autonomous Kirchhoff wave equations with strong damping and delay effects in the phase space $\mathcal{E} = C_{H^1_0(\Omega)} \times C_{L^2(\Omega)}$. The model contains a strong damping term $-\Delta \partial_t u$, a nonlocal Kirchhoff term $\Psi(\|\nabla u\|^2)$, a state-dependent delay term $\phi(t, u_t)$, and a time-dependent external force $h(x, t)$. There appear to be no results on the pullback attractors of Kirchhoff wave equations involving both strong damping and delay effects in the literature. To this end, we establish delicate uniform estimates and employ the contraction function method to prove the pullback asymptotic compactness of the associated process, which yields the existence of pullback $\mathcal{D}$-attractors.


\end{abstract}

\hspace{3mm}{\bf Keywords:} Kirchhoff wave equations; Pullback \(\mathcal{D}\)-attractors; Strong damping; Delay.

\hspace{3mm}{\bf 2020 MSC:} 35B40, 35B41, 35L05.

\section{Introduction}

\setcounter{equation}{0}

\let\thefootnote\relax\footnote{*Corresponding author: binyangdhu@163.com}
\let\thefootnote\relax\footnote{\footnotesize E-mails: yuming@dhu.edu.cn, alain.miranville@univ-lehavre.fr, hxximust2025@163.com}
In this paper, we study the following Kirchhoff wave equations with strong damping $- \Delta\partial_{t}u$ and delay $\phi(t, u_{t})$
\begin{align}
\left\{\begin{array}{ll}
\partial_{tt}u- \Delta\partial_{t} u-(1+\Psi(\|\nabla u\|^2)) \Delta u+g(u)=\phi(t, u_{t})+h(x,t) & \text { in } \Omega \times(\tau, +\infty), \\
u(x,t)=0 & \text { on } \partial \Omega\times(\tau, +\infty), \\
u(x, t)=\varphi(x, t-\tau), \,\partial_{t}u(x,t)=\partial_{t}\varphi(x, t-\tau), &\,\, x \in \Omega,\, t \in[\tau-k, \tau],
\end{array}\right.\label{1.1-7}
\end{align}
where $\Omega \subset \mathbb R^{n}\,(n\geq3)$ is a bounded domain with smooth boundary $\partial\Omega$, \(\Psi(\|\nabla u\|^2)\) is the Kirchhoff term defined as $\Psi(s)=\xi s$ for any $\xi>0$ and $s\in\mathbb R^+$, \(g(u)\) is the nonlinear term, \(\phi(t, u_t)\) is the delay function, \(h(x, t)\) is the non-autonomous external forcing term, and $k>0$ is the length of delay effects.

Firstly, we introduce some essential function spaces and their basic properties.

For any solution $u$, define the history variable
$$
u_t(\varrho)=u(t+\varrho),
\qquad
\forall\, t\in\mathbb{R}, \ \varrho\in[-k,0].
$$
Clearly, $u_t$ denotes the history segment of the solution, while $\partial_tu$ stands for the time derivative of $u$.

For a Banach space $X$, denote by $C([-k,0];X)$ the space of continuous functions from $[-k,0]$ into $X$, equipped with the norm
$$
\|\phi\|_{C([-k,0];X)}=\max_{\varrho\in[-k,0]}\|\phi(\varrho)\|_X.
$$
For simplicity, we write $C_X=C([-k,0];X)$, in particular, $C_{H_0^1}=C([-k,0];H_0^1(\Omega))$ and $C_{L^2}=C([-k,0];L^2(\Omega))$.

Define $C_{H_0^1,L^2}=C_{H_0^1(\Omega), L^2(\Omega)}=C_{H_0^1(\Omega)}\cap C^1([-k,0];L^2(\Omega))$, which is a Banach space endowed with the norm
$$
\|\phi\|_{C_{H_0^1,L^2}}^2
=\|\nabla\phi\|_{C_{L^2(\Omega)}}^2+\|\phi'\|_{C_{L^2(\Omega)}}^2,
$$
where $\phi'=\partial_t\phi$. Our phase space is defined by $\mathcal{E}=C_{H_0^1(\Omega)}\times C_{L^2(\Omega)}$, with norm
$$
\|(u_t, u_t^{\prime})\|_{\mathcal{E}}^2
=\|\nabla u_t\|_{C_{L^2(\Omega)}}^2+\|u_t'\|_{C_{L^2(\Omega)}}^2,
$$
where $u_t(\varrho)=u(t+\varrho)$ and $u_t'=\partial_t u(t+\varrho)$.

Next, we present the assumptions on the functions appearing in equation \eqref{1.1-7}.

Suppose $g(\cdot) \in C^{1}(\mathbb{R} ; \mathbb{R})$ satisfies
\begin{equation}\label{1.2-7}
\left|g^{\prime}(u)\right| \leq C\left(1+|u|^{\frac{2}{n-2}}\right),\quad \forall \,u \in \mathbb R.
\end{equation}

Define $G(u)=\int_0^u g(s) d s$, then from \eqref{1.2-7}, we deduce that there exists a constant $C>0$ such that
\begin{equation}\label{1.3-7}
|G(u)| \leq C\left(1+|u|^{\frac{2 n-2}{n-2}}\right).
\end{equation}

Moreover, assume $G(\cdot)$ satisfies the following inequalities
\begin{equation}\label{1.4-7}
\lim _{|u| \rightarrow \infty} \inf \frac{u g(u)-\chi G(u)}{u^2} \geq 0, \quad\forall \,\chi>0,
\end{equation}
\begin{equation}\label{1.5-7}
\lim _{|u| \rightarrow \infty} \inf \frac{G(u)}{u^2} \geq 0.
\end{equation}

Let \(\phi(t, u_t): \mathbb{R} \times C_{L^2(\Omega)} \to L^2(\Omega)\) with \(\phi(t, 0) = 0\), satisfying the following assumptions below:

($A_1$) For any \(v \in C_{L^2(\Omega)}\), the mapping \(t \mapsto \phi(t, v) \in L^2(\Omega)\) is measurable;

($A_2$) There exists a constant \(\widetilde{C} = \widetilde{C}(\phi) > 0\) such that for all $v,w\in C_{L^2(\Omega)}$,
\begin{equation}\label{1.6-7}
\|\phi(t, v) - \phi(t, w)\| \leq \widetilde{C}\|v - w\|_{C_{L^2(\Omega)}};
\end{equation}

($A_3$) There exist constants $\mu_0>0$ and $C_\phi>0$ such that for any $\mu\in [0,\mu_0]$ and all $v,w\in L^2([\tau-k,t];L^2(\Omega))$, the following inequality holds:
\begin{equation}\label{1.7-7}
\int_\tau^t e^{\mu s}\|\phi(s,v_s)-\phi(s,w_s)\|^2 ds
\le C_\phi^2 \int_{\tau-k}^t e^{\mu s}\|v(s)-w(s)\|^2 ds.
\end{equation}

Additionally, assume \(h(x, t) \in L_{\text{loc}}^2\left(\mathbb{R}; L^2(\Omega)\right)\), implying that it is locally square-integrable with respect to \(t \in \mathbb{R}\).

Kirchhoff-type wave equations originate from Kirchhoff \cite{k.5}, who extended the classical d'Alembert wave equation by accounting for the variation of tension caused by the deformation of elastic strings. Since then, various Kirchhoff-type models have been widely investigated from both theoretical and applied perspectives. A prototypical example is
$u_{tt}-\left(1+\delta\|\nabla u\|^2\right)\Delta u+f(u)=g(x)$ with $\delta>0$,
which has served as a fundamental model in the study of nonlinear wave propagation.

For non-autonomous Kirchhoff wave equations with strong damping, the existence and long-time behavior of pullback attractors have been investigated in \cite{wz.5, ly.5, MWX.7, lyf.5, YQAW.7} and references therein. Wang and Zhong \cite{wz.5} discussed the Kirchhoff wave equation
$u_{tt}-\Delta u_t-\left(1+\delta_1 \|\nabla u\|^2\right) \Delta u+f(u)=g(x, t)$
with $\delta_1>0$ in $H_0^1(\Omega)\times L^2(\Omega)$. By decomposing the associated process into the sum of a decaying part and a compact part, they verified the upper semicontinuity of pullback attractors with respect to the parameter $\delta_1$.
Moreover, Li and Yang \cite{ly.5} established the robustness of pullback attractors and pullback exponential attractors for a class of the Kirchhoff wave equation with strong nonlinear damping
$u_{t t}-\left(1+\delta_2\|\nabla u\|^2\right) \Delta u-\sigma(\|\nabla u\|^2) \Delta u_t+f(u)=g(x, t)$ in $(H_0^1(\Omega) \cap L^{p+1}(\Omega))\times L^2(\Omega)$ with $1 \leq p \le \frac{N+4}{(N-4)^{+}}$ and $N\geq3$, where $\delta_2 \in[0,1\rceil$ is an extensible parameter, and $\sigma \in C^1\left(\mathbb{R}^{+}\right)$ with $\sigma(s) \geq \sigma_0>0$.
By employing the measure of noncompactness method, Ma, Wang and Xie \cite{MWX.7} obtained the existence of pullback attractors for the Kirchhoff wave equation $u_{t t}-\Delta u_t-\phi\left(\|\nabla u\|^2\right) \Delta u+f(u)=h(x, t)$ in $H_0^1(\Omega) \times L^2(\Omega)$.
Later, Li, Yang and Feng \cite{lyf.5} discussed the existence and continuity of uniform attractors for the Kirchhoff wave equation $u_{t t}-\left(1+\delta_3\|\nabla u\|^2\right) \Delta u-\Delta u_t+f(u)=g(x, t)$, where the perturbation parameter $\delta_3\in[0,1]$.
In our recent work \cite{YQAW.7}, we established the existence and upper semicontinuity of pullback attractors for the Kirchhoff wave equations $\varepsilon(t) \partial_{t t} u-\left(1+\delta_4\|\nabla u\|^2\right) \Delta u-\Delta \partial_t u+\lambda u=g(u)+h(t)$ with $\delta_4, \lambda>0$ in the time-dependent space $X_t$ with the norm $\left\|\nabla u_0\right\|^2+\varepsilon(t)\|u_1\|^2$.
In addition, several results regarding attractors for delayed Kirchhoff-type equations are available in \cite{KM.7, QW.7, QYH.7, Park.7}.

However, the study of Kirchhoff-type wave equations that simultaneously involve strong damping and delay effects remains largely unexplored. To the best of our knowledge, no results on the existence or upper semicontinuity of pullback attractors for such systems have been reported.
Based on the preceding studies, we identify the key difficulties in analyzing equation \eqref{1.1-7} and present the main contributions of this paper as follows:

1. \emph{Refined Analysis of the Nonlocal Kirchhoff Term.}
Equation \eqref{1.1-7} includes the nonlocal Kirchhoff term $\Psi(\|\nabla u\|^2)$ with $\Psi(s)=\xi s$ ($\xi>0$). This nonlocal term introduces two main difficulties: first, it couples the equation globally through the gradient norm, thus standard local estimates do not directly apply; second, for large $\|\nabla u\|$, it can dominate the dynamics and complicate uniform energy bounds. To address these issues, we construct higher-order energy estimates that control the growth of $\|\nabla u\|$, and combine them with compactness arguments adapted to the nonlocal structure, enabling the analysis of the long-time behavior of solutions.

2. \emph{Comprehensive Treatment of Strong Damping and Delay.}
In addition to the difficulties posed by the nonlocal Kirchhoff term, equation \eqref{1.1-7} also contains the strong damping term $-\Delta \partial_t u$ and the state-dependent delay term $\phi(t,u_t)$, whose interaction creates further obstacles in deriving uniform estimates and establishing asymptotic compactness. In particular, the delay term depends on the history $u_t$ and is non-autonomous, which prevents the direct application of standard compactness arguments. To overcome this, under assumptions $(A_1)$-$(A_3)$, we develop suitable estimates for the delay term and incorporate it into the energy framework, allowing us to handle the coupling effect and prove the existence of pullback $\mathcal{D}$-attractors.

3. \emph{Asymptotic analysis in history space.}
The problem is formulated in the phase space $\mathcal{E}=C_{H_0^1(\Omega)}\times C_{L^2(\Omega)}$, where the presence of memory requires working with history segments $u_t$. This induces a loss of compactness under time translations, making verification of asymptotic compactness highly nontrivial. By combining uniform energy estimates with compactness arguments adapted to the history structure, we control the evolution of $u_t$ and establish the asymptotic compactness of the associated process.

The structure of this paper is as follows.
In Section 2, we shall introduce some necessary definitions and lemmas.
In Section 3, we shall systematically prove the existence of pullback \(\mathcal{D}\)-attractors for equation \eqref{1.1-7}.
First, we establish the global well-posedness of solutions to equation \eqref{1.1-7}, and then deduce the existence of the process generated by solutions. Next, we perform a priori estimates on solutions and construct an appropriate universe domain to derive the existence of the pullback \(\mathcal{D}\)-absorbing set. Subsequently, by constructing multiple test functions with suitable parameters and precisely defining the range of these parameters, we conduct refined estimates of the norm of solutions in $\mathcal E=C_{H^1_0(\Omega)}\times C_{L^2(\Omega)}$. This process eventually leads to deriving that the norm of solutions in the phase space is bounded by some constants plus a certain complex function. Then, we prove that this complex function is a contraction function, thereby establishing the asymptotic compactness of the process. Finally, by combining the above results, we arrive at the main result of this paper: the existence of pullback \(\mathcal{D}\)-attractors.

\section{Preliminaries}
In this section, we recall some basic concepts and preliminary results concerning pullback attractors; see \cite{qs.5,r,QYbook1,l.5,c.5,lw.5,Lions,ljs.5,mz.5,qy.2,yd.5,ydi.5,yw.5} for further details.
Throughout this paper, $X$ denotes a Banach space (or a closed subset of a Banach space), and $\Gamma(X)$ denotes the family of all nonempty subsets of $X$. Furthermore, let $\mathcal{D}$ be a nonempty collection of families
$$
\widetilde{D}=\{D(t):t\in\mathbb{R}\},
$$
where $D(t)\in\Gamma(X)$ for each $t\in\mathbb{R}$.

\begin{Definition}\label{def1-7}
A family of mappings
$$
{U(t,\tau):X\to X,\quad t\ge\tau,\ \tau\in\mathbb{R}}
$$
is called a \emph{process} on $X$ if $U(\tau,\tau)=I$ for every $\tau\in\mathbb{R}$ and $U(t,s)U(s,\tau)=U(t,\tau)$ for all $t\ge s\ge\tau$.
\end{Definition}

\begin{Definition}\label{def2-7}
The process $\{U(t,\tau)\}_{t\ge\tau}$ is called pullback $\mathcal D$-asymptotically compact in $X$ if, for every $t\in\mathbb R$ and every family $\widetilde D\in\mathcal D$, any sequence $\tau_n\to-\infty$ together with any choice of points $x_n\in D(\tau_n)$ generates a sequence $\{U(t,\tau_n)x_n\}_{n=1}^{\infty}$ that is relatively compact in $X$.
\end{Definition}

\begin{Definition}\label{def3-7}
A family $\widetilde B\in\mathcal D$ is said to be pullback $\mathcal D$-absorbing for the process $\{U(t,\tau)\}_{t\ge\tau}$ if, for every $\widetilde D\in\mathcal D$ and every $t\in\mathbb R$, there exists $\tau_0=\tau_0(t,\widetilde D)<t$ such that
$$
U(t,\tau)D(\tau)\subset B(t)
$$
whenever $\tau\le\tau_0(t,\widetilde D)$.
\end{Definition}

\begin{Definition}\label{def4-7}
A family
$$
\mathcal A=\{\mathcal A(t):t\in\mathbb R\}\subset\Gamma(X)
$$
is called a pullback $\mathcal D$-attractor for the process
$\{U(t,\tau)\}_{t\ge \tau}$ if the following properties hold:

\noindent (i) $\mathcal A(t)$ is compact in $X$ for every $t\in\mathbb R$;

\noindent (ii) for every $\widetilde D\in\mathcal D$ and every $t\in\mathbb R$,
$$
\lim_{\tau\to-\infty}
\operatorname{dist}_X\bigl(U(t,\tau)D(\tau),\mathcal A(t)\bigr)=0;
$$

\noindent (iii) for all $-\infty<\tau\le t<+\infty$,
$$
U(t,\tau)\mathcal A(\tau)=\mathcal A(t).
$$
\end{Definition}

\begin{Definition}\label{def5-7}
Let $B$ be a nonempty subset of $X$. A function
$$
\psi:X\times X\to\mathbb R
$$
is called a contractive function on $B\times B$ if, for every sequence $\{x_n\}_{n=1}^{\infty}\subset B$, there exists a subsequence $\{x_{n_k}\}_{k=1}^{\infty}$ such that
$$
\lim_{k\to\infty}\lim_{l\to\infty}\psi(x_{n_k},x_{n_l})=0.
$$
We denote by $\operatorname{Contr}(B)$ the family of all contractive functions on $B\times B$.
\end{Definition}

\begin{Lemma}\label{Lem1-7}
Let $\{U(t,\tau)\}_{t\ge\tau}$ be a process on $X$ and suppose that it possesses a pullback $\mathcal D$-absorbing family $\widetilde B=\{B(t):t\in\mathbb R\}$. Assume further that for every $\varepsilon>0$, there exist $T=T(t,\widetilde B,\varepsilon)>0$ and
$\psi_{t,T}(\cdot,\cdot)\in\operatorname{Contr}(B(t-T))$ such that
$$
\|U(t,t-T)x-U(t,t-T)y\|_X
\leq
\varepsilon+\psi_{t,T}(x,y)
$$
for all $x,y\in B(t-T)$.
Then $\{U(t,\tau)\}_{t\ge\tau}$ is pullback $\mathcal D$-asymptotically compact in $X$.
\end{Lemma}

\begin{Lemma}\label{LemmaLH}
Let $\{U(t,\tau)\}_{t\geq \tau}$ be a process acting on $X$. Then $\{U(t,\tau)\}_{t\geq \tau}$ admits a pullback $\mathcal{D}$-attractor in $X$ if the two conditions below are satisfied:

\noindent (i) The process $\{U(t,\tau)\}_{t\geq \tau}$ possesses a pullback $\mathcal{D}$-absorbing set $B_0\subset X$;

\noindent (ii) $\{U(t,\tau)\}_{t\geq \tau}$ is pullback $\mathcal{D}$-asymptotically compact on $\widehat{B}_0$.
\end{Lemma}

\section{Existence of Pullback \(\mathcal{D}\)-Attractors $\mathcal {A}(t)$}

In this section, we aim to establish the existence of pullback \(\mathcal{D}\)-attractors for problem \eqref{1.1-7} in $\mathcal E=C_{H^1_0(\Omega)}\times C_{L^2(\Omega)}$.
By virtue of the Faedo-Galerkin scheme, we first derive the global well-posedness of \eqref{1.1-7}, including existence, uniqueness and continuous dependence of solutions. With these results in hand, we subsequently construct the associated process.

\begin{Theorem}\label{lem3.0-7}
Suppose $(\varphi, \partial_t \varphi) \in C_{H_0^{1}(\Omega)} \times C_{L^2(\Omega)}$ is given. Then for any $t \in \mathbb{R}$ and $T > \tau$, there exists a unique weak solution $u$ satisfying
\begin{equation}\label{3.1-7}
u \in C\left([\tau - k, T]; H_0^1(\Omega)\right) \cap C^1\left([\tau - k, T]; L^2(\Omega)\right).
\end{equation}
Moreover, let $u$ and $v$ be two solutions to equation \eqref{1.1-7} corresponding to different initial data, and define $\bar{u} = u - v$. Then for any $t \geq \tau$, the following Lipschitz continuity estimate holds:
\begin{equation}\label{3.2-7}
\left\|\nabla \bar{u}_t\right\|^2_{C_{L^2(\Omega)}} + \left\|\bar{u}_t^{\prime}\right\|^2_{C_{L^2(\Omega)}} \leq e^{C(t - \tau)}\left(\left\|\nabla \bar{u}_\tau\right\|^2_{C_{L^2(\Omega)}} + \left\|\bar{u}_\tau^{\prime}\right\|^2_{C_{L^2(\Omega)}}\right),
\end{equation}
where $C = C(\widetilde{C})$ is a positive constant and $\widetilde{C}$ is the constant arising from \eqref{1.6-7}.
\end{Theorem}

\begin{Corollary}
For any weak solution $u$ to equation \eqref{1.1-7}, let $\eta_u = (u, \partial_t u)$. Then we can define the continuous process $U(t, \tau): C_{H_0^1(\Omega)} \times C_{L^2(\Omega)} \rightarrow C_{H_0^1(\Omega)} \times C_{L^2(\Omega)}$ associated with problem \eqref{1.1-7} by
\begin{equation}\label{3.3-7}
U(t, \tau)\eta_u = (u(t), \partial_t u(t)), \quad \forall\, t \geq \tau,\; \eta_u \in C_{H_0^1(\Omega)} \times C_{L^2(\Omega)}.
\end{equation}
\end{Corollary}

We now derive several a priori estimates in suitable functional spaces, which will be used to establish the existence of a pullback $\mathcal D$-absorbing set.

\begin{Lemma}\label{lem3.1-7}
Let $(\varphi, \partial_t \varphi) \in C_{H_0^{1}(\Omega)} \times C_{L^2(\Omega)}$ be given, the weak solution to \eqref{1.1-7} satisfies
\begin{equation}\label{3.4-7}
\begin{aligned}
\left\|u_{t}\right\|^2_{C_{H_0^1(\Omega)}}+\|u_t^{\prime}\|_{C_{L^2(\Omega)}}^2 & \leq C_1 e^{-\zeta(t-k-\tau)} + C_2 e^{-\zeta(t-k)} \int_\tau^t e^{\zeta s}\|h(x, s)\|^2 d s + C_3,
\end{aligned}
\end{equation}
where $t-k \geq \tau$, $\zeta=\zeta\left(\lambda_1, C_\phi\right)>0$, $C_1=C_1\left(\xi, \lambda_1, C_\phi, \varphi\right)>0$, $C_2=C_2(\lambda_1)>0$, $C_3=C_3\left(\lambda_1, C_\phi\right)>0$, and the parameter $\lambda_1$ is derived from the Poincar\'e inequality $\lambda_1\|u\|^2\leq \|\nabla u\|^2$.
\end{Lemma}
\noindent $\mathbf{Proof.}$ By the Young and Poincar\'e inequalities, there exists a constant $\alpha > 0$ such that
\begin{equation}\label{3.5-7}
(\phi(t, u_t) + h(x, t), u) \leq \frac{1}{\alpha}\|\phi(t, u_t)\|^2 + \frac{1}{\alpha}\|h(x, t)\|^2 + \frac{\alpha}{2 \lambda_1}\|\nabla u\|^2.
\end{equation}

Taking the $L^2(\Omega)$ inner product of $\varepsilon u\ (\varepsilon > 0)$ with equation \eqref{1.1-7}, and applying \eqref{3.5-7} and the Poincar\'e inequality, we derive
\begin{equation}\label{3.6-7}
\begin{aligned}
& \frac{d}{dt} \left( 2 \varepsilon (\partial_t u, u) + \varepsilon \|\nabla u\|^2 \right) + \varepsilon \left( 2 - \frac{\alpha}{\lambda_1} \right) \|\nabla u\|^2 + 2 \varepsilon \xi \|\nabla u\|^4\\
& - \frac{2 \varepsilon}{\lambda_1} \|\nabla \partial_t u\|^2 + 2 \varepsilon (g(u), u)
 \leq \frac{2 \varepsilon}{\alpha} \|\phi(t, u_t)\|^2 + \frac{2 \varepsilon}{\alpha} \|h(x, t)\|^2.
\end{aligned}
\end{equation}

Similarly as in \eqref{3.5-7}, by applying the Young and Poincar\'e inequalities, there exists a constant $\beta > 0$ such that
\begin{equation}\label{3.7-7}
(\phi(t, u_t) + h(x, t), \partial_t u) \leq \frac{1}{\beta} \|\phi(t, u_t)\|^2 + \frac{1}{\beta} \|h(x, t)\|^2 + \frac{\beta}{2 \lambda_1} \| \nabla \partial_t u \|^2.
\end{equation}

Taking the $L^2(\Omega)$ inner product of equation $\eqref{1.1-7}_1$ with $\partial_t u$, and using \eqref{3.7-7} and $G(u) = \int_0^u g(s)\, ds$, we obtain
\begin{equation}\label{3.8-7}
\begin{aligned}
&\frac{d}{dt} \left( \|\partial_t u\|^2 + \|\nabla u\|^2 + \frac{\xi}{2} \|\nabla u\|^4 + 2(G(u), 1) \right)
+ \left(2 - \frac{\beta}{\lambda_1}\right) \|\nabla \partial_t u\|^2\\
&\leq \frac{2}{\beta} \|\phi(t, u_t)\|^2 + \frac{2}{\beta} \|h(x, t)\|^2.
\end{aligned}
\end{equation}

By summing inequalities \eqref{3.6-7} and \eqref{3.8-7}, we deduce
\begin{align}\label{3.9-7}
& \frac{d}{d t}\left(\|\partial_t u\|^2+(1+\varepsilon)\|\nabla u\|^2+2 \varepsilon(\partial_t u, u)+\frac{\xi}{2}\|\nabla u\|^4+2(G(u), 1)\right)\non \\
& +\varepsilon\left(2-\frac{\alpha}{\lambda_1}\right)\|\nabla u\|^2+\left(2-\frac{\beta}{\lambda_1}-\frac{2 \varepsilon}{\lambda_1}\right)\|\nabla \partial_t u\|^2+2 \varepsilon(g(u), u) +2 \varepsilon \xi\|\nabla u\|^4\non\\
&\leq 2\left(\frac{\varepsilon}{\alpha}+\frac{1}{\beta}\right)\|\phi(t, u_t)\|^2+2\left(\frac{\varepsilon}{\alpha}+\frac{1}{\beta}\right)\|h(x, t)\|^2.
\end{align}

Using \eqref{1.4-7} and the Poincar\'e inequality, there exist constants $C_4, C_5 > 0$ such that
\begin{equation}\label{3.10-7}
2 \varepsilon(g(u), u) \geq 2 \varepsilon \chi(G(u), 1) - \frac{2 \varepsilon C_4}{\lambda_1}\|\nabla u\|^2 - 2 \varepsilon C_5.
\end{equation}

Substituting \eqref{3.10-7} into \eqref{3.9-7} yields
\begin{align}\label{3.11-7}
& \frac{d}{d t}\left(\left\|\partial_t u\right\|^2+(1+\varepsilon)\|\nabla u\|^2+2 \varepsilon\left(\partial_t u, u\right)+\frac{\xi}{2}\|\nabla u\|^4+2(G(u), 1)\right) \non\\
& +\varepsilon\left(2-\frac{\alpha}{\lambda_1}-\frac{2 C_4}{\lambda_1}\right)\|\nabla u\|^2+\left(2-\frac{\beta}{\lambda_1}-\frac{2 \varepsilon}{\lambda_1}\right)\left\|\nabla \partial_t u\right\|^2\non\\
& +2 \varepsilon \chi(G(u), 1)+2 \varepsilon \xi\|\nabla u\|^4 \non\\
&\leq 2\left(\frac{\varepsilon}{\alpha}+\frac{1}{\beta}\right)\left\|\phi\left(t, u_t\right)\right\|^2+2\left(\frac{\varepsilon}{\alpha}+\frac{1}{\beta}\right)\|h(x, t)\|^2+2 \varepsilon C_5.
\end{align}

We define $A_1(t)$ by
\begin{equation}\label{3.12-7}
A_1(t) = \|\partial_t u\|^2 + (1+\varepsilon)\|\nabla u\|^2 + 2\varepsilon(\partial_t u, u) + \frac{\xi}{2} \|\nabla u\|^4 + 2(G(u), 1).
\end{equation}

Substituting \eqref{3.12-7} into \eqref{3.11-7}, we obtain
\begin{align}\label{3.13-7}
& \frac{d}{dt} A_1(t) + \varepsilon\left(2 - \frac{\alpha}{\lambda_1} - \frac{2C_4}{\lambda_1}\right) \|\nabla u\|^2 + \left(2 - \frac{\beta}{\lambda_1} - \frac{2\varepsilon}{\lambda_1}\right) \|\nabla \partial_t u\|^2 \non\\
& + 2\varepsilon \chi (G(u), 1) + 2\varepsilon \xi \|\nabla u\|^4 \non\\
& \leq 2\left(\frac{\varepsilon}{\alpha} + \frac{1}{\beta}\right) \|\phi(t, u_t)\|^2 + 2\left(\frac{\varepsilon}{\alpha} + \frac{1}{\beta}\right) \|h(x, t)\|^2 + 2\varepsilon C_5.
\end{align}

By \eqref{1.5-7} and the Poincar\'e inequality, there exist constants $C_6, C_7 > 0$ such that
\begin{equation}\label{3.14-7}
2(G(u), 1) \geq -\frac{2C_6}{\lambda_1} \|\nabla u\|^2 - 2C_7.
\end{equation}

Substituting \eqref{3.14-7} into \eqref{3.12-7}, we deduce
\begin{equation}\label{3.15-7}
A_1(t) \geq \|\partial_t u\|^2 + \left(1 + \varepsilon - \frac{2C_6}{\lambda_1} \right) \|\nabla u\|^2 + 2\varepsilon(\partial_t u, u) + \frac{\xi}{2} \|\nabla u\|^4 - 2C_7.
\end{equation}

From \eqref{3.12-7} and \eqref{3.15-7}, if we further assume $0 < C_6 < \frac{(1+\varepsilon)\lambda_1}{2}$, then there exist constants $0 < C_8 \leq \min\left\{1,\, 1+\varepsilon - \frac{2C_6}{\lambda_1} \right\}$ and $C_9 \geq C_7$ such that
\begin{equation}\label{3.16-7}
A_1(t) \geq C_8\left( \|\nabla u\|^2 + \|\partial_t u\|^2 \right) - 2C_9.
\end{equation}

From \eqref{3.12-7}, using the Young and Poincar\'e inequalities, we obtain
\begin{equation}\label{3.17-7}
A_1(t) \leq \left(1+\varepsilon+\frac{\varepsilon}{\lambda_1}\right)\|\nabla u\|^2 + \frac{1+\varepsilon}{\lambda_1}\| \nabla \partial_tu\|^2 + \frac{\xi}{2}\| \nabla u\|^4 + 2(G(u), 1).
\end{equation}

If we further assume $\alpha, C_4$ satisfy $\alpha + 2 C_4 < 2 \lambda_1$, $\beta, \varepsilon$ satisfy $\beta + 2 \varepsilon < 2 \lambda_1$, and $\sigma$ satisfies
\begin{equation}\label{3.18-7}
0 \leq \sigma \leq \min \left\{ 4 \varepsilon, \varepsilon \chi, \frac{\varepsilon \left(2 \lambda_1 - \alpha - 2 C_4 \right)}{\lambda_1 + \varepsilon \lambda_1 + \varepsilon}, \frac{2 \lambda_1 - \beta - 2 \varepsilon}{1 + \varepsilon} \right\},
\end{equation}
then it follows that
\begin{equation}\label{3.19-7}
\frac{d}{d t} A_1(t) + \sigma A_1(t) \leq 2 \left( \frac{\varepsilon}{\alpha} + \frac{1}{\beta} \right) \|\phi(t, u_t)\|^2 + 2 \left( \frac{\varepsilon}{\alpha} + \frac{1}{\beta} \right) \|h(x, t)\|^2 + 2 \varepsilon C_5.
\end{equation}

Multiplying \eqref{3.19-7} by $e^{\sigma t}$ $(
\sigma > 0)$, and then integrating the resulting expression over $[\tau, t]$, we obtain
\begin{equation}\label{3.20-7}
\begin{aligned}
e^{\sigma t} A_1(t) & \leq e^{\sigma \tau} A_1(\tau) + 2 \left( \frac{\varepsilon}{\alpha} + \frac{1}{\beta} \right) \int_\tau^t e^{\sigma s} \left\|\phi\left(s, u_s\right)\right\|^2 \, ds\\
&+ 2 \left( \frac{\varepsilon}{\alpha} + \frac{1}{\beta} \right) \int_\tau^t e^{\sigma s} \|h(x, s)\|^2 \, ds + \frac{2 \varepsilon C_5}{\sigma} e^{\sigma t}.
\end{aligned}
\end{equation}

Noticing $\phi(t, 0) = 0$ and using \eqref{1.7-7} and the Poincar\'e inequality, we conclude
\begin{align}\label{3.21-7}
& 2\left( \frac{\varepsilon}{\alpha} + \frac{1}{\beta} \right) \int_\tau^t e^{\sigma s} \left\|\phi\left(s, u_s\right)\right\|^2 \, ds \leq 2\left( \frac{\varepsilon}{\alpha} + \frac{1}{\beta} \right) C_\phi^2 \int_{\tau-k}^t e^{\sigma s} \|u(s)\|^2 \, ds \non\\
& \leq 2\left( \frac{\varepsilon}{\alpha} + \frac{1}{\beta} \right) C_\phi^2 \int_{\tau-k}^\tau e^{\sigma s} \|u(s)\|^2 \, ds + \frac{2}{\lambda_1} \left( \frac{\varepsilon}{\alpha} + \frac{1}{\beta} \right) C_\phi^2 \int_\tau^t e^{\sigma s} \|\nabla u(s)\|^2 \, ds \non\\
& \leq \frac{2}{\sigma} \left( \frac{\varepsilon}{\alpha} + \frac{1}{\beta} \right) C_\phi^2 e^{\sigma \tau} \|\varphi\|^2 + \frac{2}{\lambda_1} \left( \frac{\varepsilon}{\alpha} + \frac{1}{\beta} \right) C_\phi^2 \int_\tau^t e^{\sigma s} \|\nabla u(s)\|^2 \, ds.
\end{align}

By substituting \eqref{3.21-7} into \eqref{3.20-7}, we derive
\begin{equation}\label{3.22-7}
\begin{aligned}
e^{\sigma t} A_1(t) & \leq e^{\sigma \tau} A_1(\tau) + \frac{2}{\sigma} \left( \frac{\varepsilon}{\alpha} + \frac{1}{\beta} \right) C_\phi^2 e^{\sigma \tau} \|\varphi\|^2 + \frac{2}{\lambda_1} \left( \frac{\varepsilon}{\alpha} + \frac{1}{\beta} \right) C_\phi^2 \int_\tau^t e^{\sigma s} \|\nabla u(s)\|^2 \, ds \\
&  + 2 \left( \frac{\varepsilon}{\alpha} + \frac{1}{\beta} \right) \int_\tau^t e^{\sigma s} \|h(x, s)\|^2 \, ds  + \frac{2 \varepsilon C_5}{\sigma} e^{\sigma t}.
\end{aligned}
\end{equation}

From \eqref{3.16-7} and \eqref{3.22-7}, we deduce that there exists a constant $C_{10} > 0$ such that
\begin{equation}\label{3.23-7}
\begin{aligned}
e^{\sigma t} A_1(t)
& \leq e^{\sigma \tau} A_1(\tau)+\frac{2}{\sigma}\left(\frac{\varepsilon}{\alpha}+\frac{1}{\beta}\right) C_\phi^2  e^{\sigma \tau}\|\varphi\|^2+\frac{2}{\lambda_1}\left(\frac{\varepsilon}{\alpha}+\frac{1}{\beta}\right) C_\phi^2 C_{10} \int_\tau^t e^{\sigma s} A_{1}(s) d s \\
& +2\left(\frac{\varepsilon}{\alpha}+\frac{1}{\beta}\right) \int_\tau^t e^{\sigma s}\|h(x, s)\|^2 d s +\left(\frac{4}{\sigma\lambda_1}\left(\frac{\varepsilon}{\alpha}+\frac{1}{\beta}\right) C_\phi^2 C_9 C_{10}+\frac{2 \varepsilon C_5}{\sigma}\right) e^{\sigma t}.
\end{aligned}
\end{equation}

Define the constants ${C_i}\,(i=11,12,13,14)$ as follows
$$
C_{11} = \frac{2}{\sigma} \left( \frac{\varepsilon}{\alpha} + \frac{1}{\beta} \right) C_\phi^2 \|\varphi\|^2, \quad
C_{12} = \frac{2}{\lambda_1} \left( \frac{\varepsilon}{\alpha} + \frac{1}{\beta} \right) C_\phi^2 C_{10},
$$
$$
C_{13} = 2 \left( \frac{\varepsilon}{\alpha} + \frac{1}{\beta} \right), \quad
C_{14} = \frac{4}{\sigma \lambda_1} \left( \frac{\varepsilon}{\alpha} + \frac{1}{\beta} \right) C_\phi^2 C_9 C_{10} + \frac{2\varepsilon C_5}{\sigma}.
$$

Then by \eqref{3.23-7}, we obtain
\begin{equation}\label{3.24-7}
\begin{aligned}
e^{\sigma t} A_1(t) & \leq e^{\sigma \tau} A_1(\tau)+C_{11} e^{\sigma \tau}+C_{12} \int_\tau^t e^{\sigma s} A_1(s) d s \\
& +C_{13} \int_{\tau}^t e^{\sigma s}\|h(x, s)\|^2 d s+C_{14} e^{\sigma t}.
\end{aligned}
\end{equation}

By the Gronwall inequality, it follows that
\begin{equation}\label{3.25-7}
\begin{aligned}
e^{\sigma t} A_1(t) & \leq e^{\sigma \tau} A_1(\tau) e^{C_{12}(t-\tau)}+C_{11} e^{\sigma \tau} e^{C_{12}(t-\tau)} \\
& +C_{13} e^{C_{12}(t-\tau)} \int_\tau^t e^{\sigma s}\left\|h\left(x, s\right)\right\|^2 d s+C_{14} e^{\sigma t} e^{C_{12}(t-\tau)}.
\end{aligned}
\end{equation}

Multiplying both sides of \eqref{3.25-7} by $e^{-\sigma t}$ yields
\begin{equation}\label{3.26-7}
\begin{aligned}
A_1(t) & \leq A_1(\tau) e^{-\left(\sigma-C_{12}\right)(t-\tau)}+C_{11} e^{-\left(\sigma-C_{12}\right)(t-\tau)} \\
& +C_{13} e^{-\left(\sigma-C_{12}\right) t} e^{-C_{12} \tau} \int_\tau^t e^{\sigma s} \| h(x, s) \|^2 d s+C_{14} e^{C_{12}(t-\tau)}.
\end{aligned}
\end{equation}

Replacing $t$ with $t + \varrho$ in \eqref{3.26-7} and applying \eqref{3.16-7}, we obtain
\begin{align}\label{3.27-7}
\|\nabla u_{t}\|_{C_{L^2(\Omega)}}^2+\|u_t^{\prime}\|_{C_{L^2(\Omega)}}^2 &\leq \frac{1}{C_8} A_1(\tau+\varrho) e^{-\left(\sigma-C_{12}\right)(t+\varrho-\tau)} +\frac{C_{11}}{C_8} e^{-\left(\sigma-C_{12}\right)(t+\varrho-\tau)} \non\\
& +\frac{C_{13}}{C_8} e^{-\left(\sigma-C_{12}\right)(t+\varrho)} e^{-C_{12}(\tau+\varrho)} \int^{t}_{\tau} e^{\sigma s}\left\|h\left(x,s\right)\right\|^2 d s \non\\
& +\frac{C_{14}}{C_8} e^{C_{12}(t+\varrho-\tau)}+\frac{2 C_9}{C_8}.
\end{align}

Letting $\zeta = \sigma - C_{12}$, which depends on $\lambda_1$ and $C_{\phi}$. Then using \eqref{3.27-7}, we obtain
\begin{align}\label{3.28-7}
\|\nabla u_t\|_{C_{L^2(\Omega)}}^2 + \|u_t^{\prime}\|_{C_{L^2(\Omega)}}^2 & \leq \frac{1}{C_8} A_1(\tau + \varrho) e^{-\zeta(t + \varrho - \tau)} + \frac{C_{11}}{C_8} e^{-\zeta(t + \varrho - \tau)} \non\\
& + \frac{C_{13}}{C_8} e^{-\zeta(t + \varrho)} e^{-C_{12}(\tau + \varrho)} \int_\tau^t e^{\sigma s} \|h(x, s)\|^2 \, ds \non\\
& + \frac{C_{14}}{C_8} e^{C_{12}(t + \varrho - \tau)} + \frac{2C_9}{C_8}.
\end{align}

Since $\varrho \in [-k, 0]$, it follows from \eqref{3.28-7} that there exist constants $C_{15} \geq \frac{1}{C_8} A_1(\tau + \varrho)$ and $C_{16} = C_{16}(\sigma, \zeta) > 0$ such that
\begin{align}\label{3.29-7}
\|\nabla u_t\|_{C_{L^2(\Omega)}}^2 + \|u_t^{\prime}\|_{C_{L^2(\Omega)}}^2 & \leq \left(\frac{C_{11}}{C_8}+C_{15}\right) e^{-\zeta(t-k-\tau)} \non\\
& +\frac{C_{13} C_{16}}{C_8} e^{-\zeta(t-k)} \int_\tau^t e^{\zeta s}\left\|h\left(x,s\right)\right\|^2 d s  \non\\
& +\frac{2 C_9+C_{14}}{C_8}.
\end{align}

Consequently, by taking $C_1 = \frac{C_{11}}{C_8} + C_{15}$, $C_2 = \frac{C_{13} C_{16}}{C_8}$  and $C_3 = \frac{2C_9 + C_{14}}{C_8}$, \eqref{3.4-7} follows directly. $\hfill$$\Box$

We now proceed to establish the existence of a pullback $\mathcal{D}$-absorbing set. To this end, we first give the definition of a tempered universe.

\begin{Definition} \label{def3.1-7}
For any parameter $\delta > 0$, let $\mathcal{D} = \mathcal{D}_\delta$ denote the class of all families of nonempty subsets $\widetilde{D} = \{D(t): t \in \mathbb{R}\} \subset {\Gamma}(C_{H_0^1(\Omega), L^2(\Omega)})$ that satisfies
\begin{equation}\label{3.30-7}
\lim_{\tau \to -\infty} \left( e^{\delta \tau} \sup_{u \in D(\tau)} \|u\|_{C_{H_0^1(\Omega), L^2(\Omega)}}^2 \right) = 0,
\end{equation}
where ${\Gamma}(C_{H_0^1(\Omega), L^2(\Omega)})$ denotes the collection of all nonempty families of subsets in $C_{H_0^1(\Omega), L^2(\Omega)}$. Then, $\mathcal{D}$ is referred to as a tempered universe in $C_{H_0^1(\Omega), L^2(\Omega)}$.
\end{Definition}

\begin{Lemma}\label{lem3.2-7}
Let $\tilde{\zeta} = \sigma - C_{12} > 0$, and suppose that
\begin{equation}\label{3.31-7}
\int_{-\infty}^t e^{\tilde{\zeta} s} \|h(x, s)\|^2 \, ds < +\infty.
\end{equation}
Then the family $\widehat{B}_0 = \{B_0(t): t \in \mathbb{R}\}$, where $B_0(t) = \mathcal{B}(0, R(t))$ is the ball in $C_{H_0^1(\Omega), L^2(\Omega)}$ centered at the origin with radius $R(t)$, and $R(t)$ satisfies
\begin{equation}\label{3.32-7}
R^2(t) = C_2 e^{-\zeta(t - k)} \int_\tau^t e^{\zeta s} \|h(x, s)\|^2 \, ds + C_{17},
\end{equation}
with $C_{17} \geq \max\{C_1, C_3\} > 0$, is a pullback $\mathcal{D}$-absorbing set for the process $\{U(t, \tau)\}_{t \geq \tau}$. Moreover, $\widehat{B}_0 \in \mathcal{D}_\delta$.
\end{Lemma}

\noindent$\mathbf{Proof.}$ By Lemma \ref{lem3.1-7}, there exists a constant $C_{17} \geq \max\{C_1, C_3\} > 0$ such that
\begin{equation}\label{3.33-7}
\|u_t\|_{C_{H_0^1(\Omega)}}^2 + \|u_t^{\prime}\|_{C_{L^2(\Omega)}}^2 \leq C_2 e^{-\zeta(t - k)} \int_{\tau}^t e^{\zeta s} \|h(x, s)\|^2 \, ds + C_{17},
\end{equation}
which implies that \eqref{3.32-7} holds.

Combining \eqref{3.31-7}$-$\eqref{3.33-7}, we deduce
$$
e^{\delta t} R^2(t) \to 0 \quad \text{as } t \to -\infty,
$$
which implies that $\widehat{B}_0 \in \mathcal{D}_\delta$.
$\hfill$$\Box$

Next, we prove the pullback $\mathcal{D}$-asymptotic compactness of the process $\{U(t, \tau)\}_{t \geq \tau}$ generated by problem \eqref{1.1-7}.

\begin{Lemma}\label{lem3.3-7}
The process $\{U(t, \tau)\}_{t \geq \tau}$ generated by problem \eqref{1.1-7} is pullback $\mathcal{D}$-asymptotically compact in $C_{H_0^1(\Omega), L^2(\Omega)}$.
\end{Lemma}
\noindent $\mathbf{Proof.}$ Suppose that $\eta_a = (a(t), \partial_t a(t))$ is a weak solution of problem \eqref{1.1-7} with initial data $(\varphi_a(x, t - \tau), \partial_t \varphi_a(x, t - \tau)) \in B_0(\tau) \times B_0(\tau)$, and $\eta_b = (b(t), \partial_t b(t))$ is a weak solution of problem \eqref{1.1-7} with initial data $(\varphi_b(x, t - \tau), \partial_t \varphi_b(x, t - \tau)) \in B_0(\tau) \times B_0(\tau)$, where $t \in [\tau - k, \tau]$.

Assume $\eta_v = (v(t), \partial_t v(t))=\eta_a-\eta_b$, then by \eqref{1.1-7}, we conclude
\begin{equation}\label{3.34-7}
\left\{\begin{array}{ll}
\partial_{tt}v- \Delta\partial_{t} v-\Delta v -\xi(\|\nabla a\|^2\Delta a-\|\nabla b\|^2\Delta b) +g(a)-g(b)\\
=\phi_a-\phi_b & \text { in } \Omega \times[\tau, +\infty), \\
v(x,t)=0 & \text { on } \partial \Omega\times[\tau, +\infty), \\
v(x, t)=\varphi_a-\varphi_b, \,\partial_{t}v(x,t)=\partial_{t}\varphi_a-\partial_{t}\varphi_b, &\,\, x \in \Omega,\, t \in[\tau-k, \tau],
\end{array}\right.
\end{equation}
where $\phi_a=\phi(t, a_{t})$, $\phi_b=\phi(t, b_{t})$, $\varphi_a=\varphi_a(x, t-\tau)$, $\varphi_b=\varphi_b(x, t-\tau)$, $\partial_{t}\varphi_a=\partial_{t}\varphi_a(x, t-\tau)$ and $\partial_{t}\varphi_b=\partial_{t}\varphi_b(x, t-\tau)$.

Taking the $L^2(\Omega)$ inner product of $\eqref{3.34-7}_1$ with $\partial_t v$, we obtain
\begin{equation}\label{3.35-7}
\begin{aligned}
& \frac{1}{2} \frac{d}{d t}\left(\|\nabla v\|^2+\|\partial_t v\|^2\right)+\|\nabla \partial_t v\|^2+\xi\left(\|\nabla a\|^2 \nabla a-\|\nabla b\|^2 \nabla b, \nabla \partial_t v\right) \\
& +(g(a)-g(b), \partial_t v)=(\phi(t, a_t)-\phi(t, b_t), \partial_t v).
\end{aligned}
\end{equation}

By the Young and Poincar\'e inequalities, there exists a constant $\theta > 0$ such that
\begin{equation}\label{3.36-7}
(\phi(t, a_t)-\phi(t, b_t), \partial_t v)
\leq \frac{1}{2\theta} \| \phi(t, a_t) - \phi(t, b_t) \|^2
+ \frac{\theta}{2\lambda_1} \| \nabla \partial_t v \|^2.
\end{equation}

Let $I_1(t)=\frac{1}{2}\left(\|\nabla v\|^2+\|\partial_t v\|^2\right)$, and substituting \eqref{3.36-7} into \eqref{3.35-7}, we deduce
\begin{equation}\label{3.37-7}
\begin{aligned}
&\frac{d}{dt} I_1(t) + \left(1 - \frac{\theta}{2\lambda_1}\right) \|\nabla \partial_t v\|^2 + \xi\left(\|\nabla a\|^2 \nabla a - \|\nabla b\|^2 \nabla b, \nabla \partial_t v\right)\\
& + (g(a) - g(b), \partial_t v)
\leq \frac{1}{2\theta} \left\| \phi(t, a_t) - \phi(t, b_t) \right\|^2.
\end{aligned}
\end{equation}

Multiplying \eqref{3.37-7} by $e^{\sigma_1 t}\,(\sigma_1>0)$ and integrating the resulting inequality over the interval $[s, t]$, we arrive at
\begin{align}\label{3.38-7}
& e^{\sigma_1 t} I_1(t)-e^{\sigma_1 s} I_1(s)-\sigma_1 \int_s^t e^{\sigma_1 r} I_1(r) d r+\left(1-\frac{\theta}{2 \lambda_1}\right) \int_s^t e^{\sigma_1 r}\left\|\nabla \partial_t v(r)\right\|^2 d r \non\\
& +\xi \int_s^t e^{\sigma_1 r}\left(\|\nabla a(r)\|^2 \nabla a(r)-\|\nabla b(r)\|^2 \nabla b(r), \nabla \partial_t v(r)\right) d r \non\\
& +\int_s^t e^{\sigma_1 r}\left(g(a(r))-g(b(r)), \partial_t v(r)\right) d r \non\\
& \leq \frac{1}{2 \theta} \int_s^t e^{\sigma_1 r}\left\|\phi\left(r, a_r\right)-\phi\left(r, b_r\right)\right\|^2 d r .
\end{align}

Integrating \eqref{3.38-7} with respect to $s$ over $[\tau, t]$, we conclude
\begin{align}\label{3.39-7}
& (t-\tau) e^{\sigma_1 t} I_1(t)-\int_\tau^t e^{\sigma_1 s} I_1(s) d s +\left(1-\frac{\theta}{2 \lambda_1}\right) \int_\tau^t \int_s^t e^{\sigma_1 r}\|\nabla \partial_t v(r)\|^2 d r d s  \non\\
& +\xi \int_\tau^{t}\int_s^t e^{\sigma_1 r}\left(\|\nabla a(r)\|^2\nabla  a(r)-\|\nabla b(r)\|^2 \nabla b(r), \nabla \partial_t v(r) )d r d s \right. \non \\
& +\int_\tau^t \int_s^t e^{\sigma_1 r}(g(a(r))-g(b(r)), \partial_t v(r)) d r d s  \non\\
& \leq \sigma_1 \int_\tau^t \int_s^t e^{\sigma_1 r} I_1(r) d r d s+\frac{1}{2 \theta} \int_\tau^t \int_s^t e^{\sigma_1 r}\left\|\phi\left(r, a_r\right)-\phi\left(r, b_r\right)\right\|^2 d r d s.
\end{align}

Taking the inner product of $\eqref{3.34-7}_1$ with $v$ in $L^2(\Omega)$, we obtain
\begin{equation}\label{3.40-7}
\begin{aligned}
& \frac{d}{dt} \left( \frac{1}{2} \|\nabla v\|^2 + (\partial_t v, v) \right) - \|\partial_t v\|^2 + \|\nabla v\|^2 + \xi\left( \|\nabla a\|^2 \nabla a - \|\nabla b\|^2 \nabla b, \nabla v \right) \\
& + (g(a) - g(b), v) = (\phi(t, a_t) - \phi(t, b_t), v).
\end{aligned}
\end{equation}

Using the Young and Poincar\'e inequalities, we derive
\begin{equation}\label{3.41-7}
\left(\phi\left(t, a_t\right)-\phi\left(t, b_t\right), v\right) \leq \frac{1}{2 \lambda_1}\left\|\phi\left(t, a_t\right)-\phi\left(t, b_t\right)\right\|^2+\frac{1}{2}\|\nabla v\|^2.
\end{equation}

Substituting \eqref{3.41-7} into \eqref{3.40-7}, we arrive at
\begin{equation}\label{3.42-7}
\begin{aligned}
& \frac{d}{d t}\left(\frac{1}{2}\|\nabla v\|^2+(\partial_t v, v)\right)-\left\|\partial_t v\right\|^2+\frac{1}{2}\|\nabla v\|^2+\xi\left(\|\nabla a\|^2 \nabla a-\|\nabla b\|^2 \nabla b, \nabla v\right)\\
& +(g(a)-g(b), v) \leq \frac{1}{2 \lambda_1}\|\phi(t, a_t)-\phi(t, b_t)\|^2.
\end{aligned}
\end{equation}

Multiplying both sides of \eqref{3.42-7} by $e^{\sigma_1 t}$ yields
\begin{align}\label{3.43-7}
& \frac{d}{d t}\left(\frac{1}{2} e^{\sigma_1 t}\|\nabla v\|^2+e^{\sigma_1 t}(\partial_t v, v)\right)+\frac{1-\sigma_1}{2} e^{\sigma_1 t}\|\nabla v\|^2-\sigma_1 e^{\sigma_1 t}(\partial_t v, v) \non\\
& -e^{\sigma_1 t}\|\partial_t v\|^2+\xi e^{\sigma_1 t}\left(\|\nabla a\|^2 \nabla a-\|\nabla b\|^2 \nabla b, \nabla v\right)  +e^{\sigma_1 t}(g(a)-g(b), v) \non\\
& \leq \frac{1}{2 \lambda_1} e^{\sigma_1 t}\left\|\phi\left(t, a_t\right)-\phi\left(t, b_t\right)\right\|^2.
\end{align}

Integrating \eqref{3.43-7} over $[s, t]$ with respect to $t$ yields
\begin{align}\label{3.44-7}
& \frac{1}{2} e^{\sigma_1 t}\|\nabla v(t)\|^2+e^{\sigma_1 t}(\partial_tv(t), v(t))+\frac{1-\sigma_1}{2} \int_s^t e^{\sigma_1 r}\|\nabla v(r)\|^2 d r \non\\
& -\sigma_1 \int_s^t e^{\sigma_1 r}(\partial_tv(r), v(r)) d r-\int_s^t e^{\sigma_1 r}\|\partial_tv(r)\|^2 d r \non\\
& +\xi\int_s^t e^{\sigma_1 r}\left(\|\nabla a(r)\|^2 \nabla a(r)-\|\nabla b(r)\|^2 \nabla b(r), \nabla v(r)\right) d r \non\\
& +\int_s^t e^{\sigma_1 r}(g(a(r))-g(b(r)), v(r)) d r \non\\
& \leq \frac{1}{2} e^{\sigma_1 s}\|\nabla v(s)\|^2+e^{\sigma_1s}(\partial_tv(s), v(s)) +\frac{1}{2 \lambda_1} \int_s^t e^{\sigma_1 r}\left\|\phi\left(r, a_r\right)-\phi\left(r, b_r\right)\right\|^2 d r.
\end{align}

Integrating \eqref{3.44-7} over $[\tau, t]$ with respect to $s$ leads to
\begin{align}\label{3.45-7}
& \frac{t-\tau}{2} e^{\sigma_1 t}\|\nabla v(t)\|^2+(t-\tau) e^{\sigma_1 t}(\partial_t v(t), v(t))+\frac{1-\sigma_1}{2}{ \int_\tau^t \int_s^t e^{\sigma_1 r}\|\nabla v(r)\|^2 d r d s} \non\\
& -\sigma_1 \int_\tau^t \int_s^t e^{\sigma_1 r}(\partial_t v(r), v(r)) d r d s-\int_\tau^t \int_s^t e^{\sigma_1 r}\|\partial_t v(r)\|^2 d r d s \non\\
& +\xi \int_\tau^t \int_s^t e^{\sigma_1 r}\left(\|\nabla a(r)\|^2 \nabla a(r)-\|\nabla b(r)\|^2 \nabla b(r), \nabla v(r)\right) d r d s \non\\
& +\int_t^t \int_s^t e^{\sigma_1 r}(g(a(r))-g(b(r)), v(r)) d r d s \non\\
& \leq \frac{1}{2} \int_\tau^t e^{\sigma_1 s}\|\nabla v(s)\|^2 d s+\int_\tau^t e^{\sigma_1 s}(\partial_t v(s), v(s)) d s \non\\
& +\frac{1}{2 \lambda_1} \int_\tau^t \int_s^t e^{\sigma_1 r}\left\|\phi\left(r, a_r\right)-\phi\left(r, b_r\right)\right\|^2 d r d s.
\end{align}

Suppose $\sigma_1 \in (0,1)$, then multiplying both sides of \eqref{3.45-7} by $\sigma_2=\frac{\sigma_1}{1-\sigma_1}$, we obtain
\begin{align}\label{3.46-7}
& \frac{(t-\tau){\sigma_2}}{2} e^{\sigma_1 t}\|\nabla v(t)\|^2+\sigma_2(t-\tau) e^{\sigma_1 t}\left(\partial_{t}v(t), v(t)\right)+\sigma_1 \int_\tau^t \int_s^t e^{\sigma_1 r} I_1(r) d r d s \non\\
& -\sigma_1 \sigma_2 \int_\tau^t \int_s^t e^{\sigma_1 r}(\partial_t v(r), v(r)) d r d s -{\frac{\sigma_1(3-\sigma_1)}{2(1-\sigma_1)} \int_\tau^t \int_s^t e^{\sigma_1 r}\left\|\partial_t v(r) \|^2 d r d s\right.} \non\\
& +\sigma_2 \xi\int_\tau^t \int_s^t e^{\sigma_1 r}\left(\|\nabla a(r)\|^2 \nabla a(r)-\|\nabla b(r)\|^2 \nabla b(r), \nabla v(r)\right) d r d s \non\\
& +\sigma_2 \int_\tau^t \int_s^t e^{\sigma_1 r}(g(a(r))-g(b(r)), v(r)) d r d s \non\\
& \leq \frac{\sigma_2}{2} \int_\tau^t e^{\sigma_1 s}\|\nabla v(s)\|^2 d s+\sigma_2 \int_\tau^t e^{\sigma_1 s}(\partial_t v(s), v(s)) d s \non\\
& +\frac{\sigma_2}{2 \lambda_1} \int_\tau^t \int_s^t e^{\sigma_1 r}\left\|\phi\left(r, a_r\right)-\phi(r, b_r)\right\|^2 d r d s.
\end{align}

Then by the Poincar\'e inequality, it follows that
\begin{align}\label{3.47-7}
& \frac{(t-\tau) \sigma_2}{2} e^{\sigma_1 t}\|\nabla v(t)\|^2+\sigma_2(t-\tau) e^{\sigma_1 t}\left(\partial_t v(t), v(t)\right)+\sigma_1 \int_\tau^t \int_s^t e^{\sigma_1 r} I_1(r) d r d s \non\\
& -\sigma_1 \sigma_2 \int_\tau^t \int_s^t e^{\sigma_1 r}\left(\partial_t v(t), v(t)\right) d r d s-\frac{\sigma_1\left(3-\sigma_1\right)}{2\left(1-\sigma_1\right) \lambda_1} \int_\tau^t \int_s^t e^{\sigma_1 r}\left\|\nabla \partial_t v(t)\right\|^2 d r d s \non\\
& +\sigma_2 \xi \int_\tau^t \int_s^t e^{\sigma_1 r}\left(\|\nabla a(r)\|^2 \nabla a(r)-\|\nabla b(r)\|^2 \nabla b(r), \nabla v(r)\right) d r d s \non\\
& +\sigma_2 \int_\tau^t \int_s^t e^{\sigma_1 r}(g(a(r))-g(b(r)), v(r)) d r d s \non\\
& \leq \frac{\sigma_2}{2} \int_\tau^t e^{\sigma_1 s}\|\nabla v(s)\|^2 d s+\sigma_2 \int_\tau^t e^{\sigma_1 s}\left(\partial_t v(s), v(s)\right) d s \non\\
& +\frac{\sigma_2}{2 \lambda_1} \int_\tau^t \int_s^t e^{\sigma_1 r}\left\|\phi\left(r, a_r\right)-\phi\left(r, b_r\right)\right\|^2 d r d s.
\end{align}

Summing \eqref{3.39-7} and \eqref{3.47-7}, we conclude
\begin{align}\label{3.48-7}
&(t-\tau) e^{\sigma_1 t} I_1(t)-\int_\tau^t e^{\sigma_1 s} I_1(s) d s \non\\
&+\left(1-\frac{\sigma_1\left(3-\sigma_1\right)}{2\left(1-\sigma_1\right) \lambda_1}-\frac{\theta}{2 \lambda_1}\right) \int_\tau^t \int_s^t e^{\sigma_1 r} \| \nabla \partial_t v(r)\|^2 d r d s \non\\
& +\sigma_2 \xi \int_\tau^t \int_s^t e^{\sigma_1 r}\left(\|\nabla a(r)\|^2 \nabla a(r)-\|\nabla b(r)\|^2 \nabla b(r), \nabla v(r)\right) d r d s \non\\
& +\sigma_2 \int_\tau^t \int_s^t e^{\sigma_1 r}(g(a(r))-g(b(r)), v(r)) d r d s \non\\
& +\xi \int_\tau^t \int_s^t e^{\sigma_1 r}\left(\|\nabla a(r)\|^2 \nabla a(r)-\|\nabla b(r)\|^2 \nabla b(r), \nabla \partial_t v(r)\right) d r d s \non\\
& +\int_\tau^t \int_s^t e^{\sigma_1 r}(g(a(r))-g(b(r)), \partial_t v(r)) d r d s \non\\
& \leq-\frac{(t-\tau) {\sigma_2}}{2} e^{\sigma_1 t}\|\nabla v(t)\|^2-\sigma_2(t-\tau) e^{\sigma_1 t}(\partial_t v(t), v(t)) \non\\
& +\frac{\sigma_2}{2} \int_\tau^t e^{\sigma_1 s}\|\nabla v(s)\|^2 d s+\sigma_2 \int_\tau^t e^{\sigma_1 s}(\partial_t v(s), v(s)) d s \non\\
& +\sigma_1 \sigma_2 \int_\tau^t \int_s^t e^{\sigma_1 r}(\partial_t v(r), v(r)) d r d s \non\\
&+\left(\frac{1}{2 \theta}+\frac{\sigma_2}{2 \lambda_1}\right)\int_\tau^t\left.\int_{s}^{t} e^{\sigma_1 r}\left\|\phi\left(r, a_r\right)-\phi\left(r, b_r\right)\right\|^2 d r d s.\right.
\end{align}

Assuming that
$$
\theta \in\left(0, \frac{\sigma_1^2-\left(2 \lambda_1+3\right) \sigma_1+2 \lambda_1}{1-\sigma_1}\right)
$$
and
$$
0<\sigma_1<\min \left\{1, \frac{2 \lambda_1+3-\sqrt{4 \lambda_1^2+4 \lambda_1+9}}{2}\right\} \quad \text { or } \quad \sigma_1>\frac{2 \lambda_1+3+\sqrt{4 \lambda_1^2+4 \lambda_1+9}}{2},
$$
we deduce
\begin{align} \label{3.49-7}
\left(1-\frac{\sigma_1\left(3-\sigma_1\right)}{2\left(1-\sigma_1\right) \lambda_1}-\frac{\theta}{2 \lambda_1}\right) \int_\tau^t \int_s^t e^{\sigma_1 r}\left\|\nabla \partial_t v(r)\right\|^2 d r d s \geq 0.
\end{align}

From \eqref{3.48-7} and \eqref{3.49-7}, we conclude
\begin{align} \label{3.50-7}
& (t-\tau) e^{\sigma_1 t} I_1(t)-\int_\tau^t e^{\sigma_1 s} I_1(s) d s \non\\
& +\sigma_2\xi \int_\tau^t \int_s^t e^{\sigma_1r}\left(\|\nabla a(r)\|^2 \nabla a(r)-\|\nabla b(r)\|^2 \nabla b(r), \nabla v(r)\right) d r d s\non\\
& +\sigma_2 \int_\tau^t \int_s^t e^{\sigma_1r}(g(a(r))-g(b(r)), v(r)) d r d s \non\\
& +\xi\int_\tau^t \int_s^t e^{\sigma_1r}\left(\|\nabla a(r)\|^2 \nabla a(r)-\|\nabla b(r)\|^2 \nabla b(r), \nabla \partial_t v(r)\right) d r d s \non\\ & +\int_\tau^t \int_s^t e^{\sigma_1r}(g(a(r))-g(b(r)), \partial_t v(r)) d r d s\non\\
&\leq-\frac{(t-\tau){\sigma_2}}{2} e^{\sigma_1 t}\|\nabla v(t)\|^2-\sigma_2(t-\tau) e^{\sigma_1 t}(\partial_t v(t), v(t)) \non\\ & +\frac{\sigma_2}{2} \int_\tau^t e^{\sigma_1 s}\|\nabla v(s)\|^2 d s+\sigma_2 \int_\tau^t e^{\sigma_1 s}(\partial_t v(s), v(s)) d s\non\\
& +\sigma_1 \sigma_2 \int_\tau^t \int_s^t e^{\sigma_1 r}(\partial_t v(r), v(r)) d r d s \non\\
& +\left(\frac{1}{2 \theta}+\frac{\sigma_2}{2 \lambda_1}\right) \int_\tau^t \int_s^t e^{\sigma_1 r}\left\|\phi\left(r, a_r\right)-\phi\left(r, b_r\right)\right\|^2 d r d s.
\end{align}

Integrating \eqref{3.50-7} over $[\tau,t]$, we deduce
\begin{align}\label{3.51-7}
& \frac{1}{2} e^{\sigma_1 t}\|\nabla v(t)\|^2+e^{\sigma_1 t}(\partial_tv(t), v(t))+\frac{1-\sigma_1}{2} \int_\tau^t e^{\sigma_1 s}\|\nabla v(s)\|^2 d s\non\\
&-\sigma_1 \int_\tau^t e^{\sigma_1 s}(\partial_t v(s), v(s)) d s
-\int_\tau^t e^{\sigma_1 s}\|\partial_t v(s)\|^2 d s \non\\
&+\xi\int_\tau^t e^{\sigma_1 s}(\|\nabla a(s)\|^2 \nabla a(s)-\|\nabla b(s)\|^2 \nabla b(s), \nabla v(s)) d s \non\\
& +\int_\tau^t e^{\sigma_1 s}(g(a(s))-g(b(s)), v(s)) d s \non\\
& \leq \frac{1}{2} e^{\sigma_1 \tau}\|\nabla v(\tau)\|^2+e^{\sigma_1 \tau}(\partial_t v(\tau), v(\tau))\non\\
&+\frac{1}{2 \lambda_1} \int_\tau^t e^{\sigma_1 s}\left\|\phi\left(s, a_s\right)-\phi(s, b_s)\right\|^2 d s.
\end{align}

Multiplying \eqref{3.51-7} by $\sigma_3=\frac{2}{1-\sigma_1}$ and then adding $\displaystyle \int_{\tau}^{t} e^{\sigma_1 s}\|\partial_t v(s)\|^{2}ds$ to both sides of the resulting inequality, we obtain
\begin{align}\label{3.52-7}
&\frac{\sigma_3}{2} e^{\sigma_1 t}\|\nabla v(t)\|^2+\sigma_3 e^{\sigma_1 t}(\partial_t v(t), v(t))\non\\
&+2 \int_\tau^t e^{\sigma_1 s} I_1(s) d s -\sigma_1 \sigma_3 \int_\tau^t e^{\sigma_1 s}(\partial_t v(s), v(s)) d s \non\\
& +\sigma_3\xi\int_\tau^t e^{\sigma_1 s}\left(\|\nabla a(s)\|^2 \nabla a(s)-\|\nabla b(s)\|^2 \nabla b(s), \nabla v(s)\right) d s \non\\
& +\sigma_3 \int_\tau^t e^{\sigma_1 s}(g(a(s))-g(b(s)), v(s)) d s \non\\
& \leq \frac{\sigma_3}{2} e^{\sigma_1 \tau}\|\nabla v(\tau)\|^2+\sigma_3 e^{\sigma_1 \tau}(\partial_t v(\tau), v(\tau))+\left(1+\sigma_3\right) \int_\tau^t e^{\sigma_1 s}\|\partial_t v(s)\|^2 d s\non\\
&+\frac{\sigma_3}{2 \lambda_1} \int_\tau^t e^{\sigma_1 s}\left\|\phi\left(s, a_s\right)-\phi(s, b_s)\right\|^2 d s.
\end{align}

Combining \eqref{3.50-7} and \eqref{3.52-7}, we obtain
\begin{align}\label{3.53-7}
&(t-\tau) e^{\sigma_1 t} I_1(t)+\int_\tau^t e^{\sigma_1 s} I_1(s) d s+\frac{(t-\tau) \sigma_2+\sigma_3}{2} e^{\sigma_1 t}\|\nabla v(t)\|^2 \non\\
&+\sigma_2 \xi \int_\tau^t \int_s^t e^{\sigma_1 r}\left(\|\nabla a(r)\|^2 \nabla a(r)-\|\nabla b(r)\|^2 \nabla b(r), \nabla v(r)\right) d r d s \non\\
&+\xi\int_\tau^t \int_s^t e^{\sigma_1 r}\left(\|\nabla a(r)\|^2 \nabla a(r)-\|\nabla b(r)\|^2 \nabla b(r), \nabla \partial_t v(r)\right) d r d s \non\\
&+\sigma_3\xi \int_\tau^t e^{\sigma_1 s}(\|\nabla a(s)\|^2 \nabla a(s)-\|\nabla b(s)\|^2 \nabla b(s), \nabla v(s)) d s \non\\
& +\sigma_2 \int_\tau^t \int_s^t e^{\sigma_1 r}(g(a(r))-g(b(r)), v(r)) d r d s +\int_\tau^t \int_s^t e^{\sigma_1 r}(g(a(r))-g(b(r)), \partial_t v(r)) d r d s \non\\
&+\sigma_3 \int_\tau^t e^{\sigma_1 s}(g(a(s))-g(b(s)), v(s)) d s\non\\
&\leq \frac{\sigma_3}{2} e^{\sigma_1 \tau}\|\nabla v(\tau)\|^2+\sigma_3 e^{\sigma_1 \tau}(\partial_t v(\tau), v(\tau)) -
\left(\sigma_2(t-\tau)+\sigma_3\right) e^{\sigma_1 t}(\partial_t v(t), v(t)) \non\\
&+\left(\sigma_1 \sigma_3+\sigma_2\right) \int_\tau^t e^{\sigma_1 s}(\partial_t v(s), v(s)) d s
+\sigma_1 \sigma_2 \int_\tau^t \int_s^t e^{\sigma_1 r}(\partial_t v(r), v(r)) d r d s \non\\
&+\frac{\sigma_2}{2} \int_\tau^t e^{\sigma_1 s}\|\nabla v(s)\|^2 d s+\left(1+\sigma_3\right) \int_\tau^t e^{\sigma_1 s}\|\partial_t v(s)\|^2 d s \non\\
&+\left(\frac{1}{2 \theta}+\frac{\sigma_2}{2 \lambda_1}\right) \int_\tau^t \int_s^t e^{\sigma_1 r}\left\|\phi\left(r, a_r\right)-\phi\left(r, b_r\right)\right\|^2 d r d s\non\\
&+\frac{\sigma_3}{2 \lambda_1} \int_\tau^t e^{\sigma_1 s}\left\|\phi\left(s, a_s\right)-\phi\left(s, b_s\right)\right\|^2 d s.
\end{align}

Multiplying \eqref{3.53-7} by $e^{\sigma_1 t}$, integrating the resulting expression over $[\tau, t]$ and applying the Poincar\'e inequality, we obtain
\begin{align}\label{3.54-7}
& e^{\sigma_1t} I_1(t)+\frac{2 \lambda_1-\theta}{2} \int_\tau^t e^{\sigma_1 s}\|\partial_t v(s)\|^2 d s \non\\ & +\xi\int_\tau^t e^{\sigma_1s}\left(\|\nabla a(s)\|^2 \nabla a(s)-\|\nabla b(s) \|^2 \nabla b(s), \nabla \partial_t v(s)\right) d s\non\\
& +\int_\tau^t e^{\sigma_1 s}(g(a(s))-g(b(s)), \partial_t v(s)) d s \non\\
& \leq e^{\sigma_1 \tau} I_1(\tau)+\sigma_1 \int_\tau^t e^{\sigma_1 s} I_1(s) d s \non\\
& +\frac{1}{2 \theta} \int_\tau^t e^{\sigma_1 s}\left\|\phi\left(s, a_s\right)-\phi\left(s, b_s\right)\right\|^2 d s.
\end{align}

Let $\theta_1=\frac{2\left(1+\sigma_3\right)}{2 \lambda_1-\theta}$ and multiplying \eqref{3.54-7} by $\theta_1$, we deduce
\begin{align}\label{3.55-7}
& \left(1+\sigma_3\right) \int_\tau^t e^{\sigma_1 s}\| \partial_ tv(s) \|^2 d s \non\\
& \leq\theta_1 e^{\sigma_1 \tau} I_1(\tau)+\sigma_1 \int_\tau^t e^{\sigma_1 s} I_1(s) d s+\frac{\theta_1}{2 \theta} \int_\tau^t e^{\sigma_1 s}\left\|\phi\left(s, a_s\right)-\phi\left(s, b_s\right)\right\|^2 d s\non\\
& -\theta_1 \xi\int_\tau^t e^{\sigma_1 s}\left(\|\nabla a(s)\|^2 \nabla a(s)-\|\nabla b(s)\|^2 \nabla b(s), \nabla \partial_t v(s)\right) d s\non \\
& -\theta_1 \int_\tau^t e^{\sigma_1 s}(g(a(s))-g(b(s)), \partial_t v(s)) d s.
\end{align}

Substituting \eqref{3.55-7} into \eqref{3.53-7} yields
\begin{align}\label{3.56-7}
&(t-\tau) e^{\sigma_1 t} I_1(t)+\left(1-\sigma_1\right) \int_\tau^t e^{\sigma_1 s} I_1(s) d s+\frac{(t-\tau) \sigma_2+\sigma_3}{2} e^{\sigma_1 t}\|\nabla v(t)\|^2\non\\
&+\sigma_2 \xi \int_\tau^t \int_s^t e^{\sigma_1 r}\left(\|\nabla a(r)\|^2 \nabla a(r)-\|\nabla b(r)\|^2 \nabla b(r), \nabla v(r)\right) d r d s \non\\
& +\xi\int_\tau^t \int_s^t e^{\sigma_1 r}\left(\|\nabla a(r)\|^2 \nabla a(r)-\|\nabla b(r)\|^2 \nabla b(r), \nabla \partial_t v(r)\right) d r d s \non\\
& \left.+\sigma_3\xi\int_\tau^t e^{\sigma_1 s}(\|\nabla a(s)\|^2 \nabla a(s)-\|\nabla b(s)\|^2 \nabla b(s), \nabla v(s)) d s\right.\non\\
& +\theta_1\xi \int_\tau^t e^{\sigma_1 s}\left(\|\nabla a(s)\|^2 \nabla a(s)-\|\nabla b(s)\|^2 \nabla b(s), \nabla \partial_t v(s)\right) d s \non\\
& +\sigma_2 \int_\tau^t \int_s^t e^{\sigma_1 r}(g(a(r))-g(b(r)), v(r)) d r d s
 +\int_\tau^t \int_s^t e^{\sigma_1 r}(g(a(r))-g(b(r)), \partial_t v(r)) d r d s \non\\
& +\sigma_3 \int_\tau^t e^{\sigma_1 s}(g(a(s))-g(b(s)), v(s)) d s
+\theta_1 \int_\tau^t e^{\sigma_1 s}(g(a (s))-g(b(s)), \partial_t v(s)) d s\non\\
& \leq  \frac{\sigma_3}{2} e^{\sigma_1 \tau}\|\nabla v(\tau)\|^2+\sigma_3 e^{\sigma_1 \tau}(\partial_t v(\tau), v(\tau))+\theta_1 e^{\sigma_1 \tau} I_1(\tau) \non\\
& -\left(\sigma_2(t-\tau)+\sigma_3\right) e^{\sigma_1 t}(\partial_t v(t), v(t))+\left(\sigma_1 \sigma_3+\sigma_2\right) \int_\tau^t e^{\sigma_1 s}(\partial_t v(s), v(s)) d s \non\\ & +\sigma_1 \sigma_2 \int_\tau^t \int_s^t e^{\sigma_1 r}(\partial_t v(r), v(r)) d r d s
+\frac{\sigma_2}{2} \int_\tau^t e^{\sigma_1s}\|\nabla v(s)\|^2 d s\non\\
&+\left(\frac{1}{2 \theta}+\frac{\sigma_2}{2 \lambda_1}\right) \int_\tau^t \int_s^t e^{\sigma_1r}\left\|\phi\left(r, a_r\right)-\phi\left(r, b_r\right)\right\|^2 d r d s \non\\
& +\left(\frac{\theta_1}{2 \theta}+\frac{\sigma_3}{2 \lambda_1}\right) \int_\tau^t e^{\sigma_1 s}\left\|\phi\left(s, a_s\right)-\phi\left(s, b_s\right)\right\|^2 d s.
\end{align}

Moreover, it follows from Lemma \ref{lem3.2-7} that there exists a constant $C_{18}=C_{18}\left(\sigma_2\right)>0$ such that
\begin{align}\label{3.57-7}
& I_1(t)+\frac{1-\sigma_1}{t-\tau} e^{-\sigma_1 t} \int_\tau^t e^{\sigma_1 s} I_1(s) d s \non\\
& +\frac{\sigma_2 \xi}{t-\tau} e^{-\sigma_1 t} \int_\tau^t \int_s^t e^{\sigma_1 r}\left(\|\nabla a(r)\|^2 \nabla a(r)-\|\nabla b(r)\|^2 \nabla b(r), \nabla v(r)\right) d r d s \non\\
& +\frac{\xi}{t-\tau} e^{-\sigma_1 t} \int_\tau^t \int_s^t e^{\sigma_1 r}\left(\|\nabla a(r)\|^2 \nabla a(r)-\|\nabla b(r)\|^2 \nabla b(r), \nabla \partial_t v(r)\right) d r d s\non \\
& +\frac{\sigma_3\xi}{t-\tau} e^{-\sigma_1 t} \int_\tau^t e^{\sigma_1 s}\left(\|\nabla a(s)\|^2 \nabla a(s)-\|\nabla b(s)\|^2 \nabla b(s), \nabla v(s)\right) d s \non\\
& +\frac{\theta_1 \xi}{t-\tau} e^{-\sigma_1 t} \int_\tau^t e^{\sigma_1 s}\left(\|\nabla a(s)\|^2 \nabla a(s)-\|\nabla b(s)\|^2 \nabla b(s), \nabla \partial_t v(s)\right) d s \non\\
& +\frac{\sigma_2}{t-\tau} e^{-\sigma_1 t} \int_\tau^t \int_s^t e^{\sigma_1r}(g(a(r))-g(b(r)), v(r)) d r d s \non\\
& +\frac{1}{t-\tau} e^{-\sigma_1 t} \int_\tau^t \int_s^t e^{\sigma_1r}(g(a(r))-g(b(r)), \partial_t v(r)) d r d s \non\\ & +\frac{\sigma_3}{t-\tau} e^{-\sigma_1 t} \int_\tau^t e^{\sigma_1 s}(g(a(s))-g(b(s)), v(s)) d s \non\\ & +\frac{\theta_1}{t-\tau} e^{-\sigma_1 t} \int_\tau^t e^{\sigma_1 s}(g(a(s))-g(b(s)), \partial_t v(s)) d s\non\\
&\leq  \frac{\sigma_3}{2(t-\tau)} e^{-\sigma_1(t-\tau)}\|\nabla v(\tau)\|^2+\frac{\sigma_3}{t-\tau} e^{-\sigma_1(t-\tau)}(\partial_t v(\tau), v(\tau))+\frac{\theta_1}{t-\tau} e^{-\sigma_1(t-\tau)} I_1(\tau) \non\\
& -\left(\sigma_2+\frac{\sigma_3}{t-\tau}\right)(\partial_t v(t), v(t))+\frac{\sigma_1 \sigma_2+\sigma_2}{t-\tau} e^{-\sigma_1 t} \int_\tau^t e^{\sigma_1 s}(\partial_t v(s), v(s)) d s \non\\
& +\frac{\sigma_1 \sigma_2}{t-\tau} e^{-\sigma_1 t} \int_\tau^t \int_s^t e^{\sigma_1 r}(\partial_t v(r), v(r)) d r d s \non\\
& +\frac{1}{t-\tau}\left(\frac{1}{2 \theta}+\frac{\sigma_2}{2 \lambda_1}\right) e^{-\sigma_1 t} \int_\tau^t \int_s^t e^{\sigma_1 r}\left\|\phi\left(r, a_r\right)-\phi\left(r, b_r\right)\right\|^2 d r d s \non\\
& +\frac{1}{t -\tau}\left(\frac{\theta_1}{2 \theta}+\frac{\sigma_3}{2 \lambda_1}\right) e^{-\sigma_1 t} \int_\tau^t e^{\sigma_1 s}\left\|\phi\left(s, a_s\right)-\phi\left(s, b_s\right)\right\|^2 ds  \non\\
& +\frac{C_{18}}{t-\tau} e^{-\sigma_1 t}.
\end{align}

Then using \eqref{1.7-7}, we conclude
\begin{align}\label{3.58-7}
&e^{-\sigma_1t} \int_\tau^t e^{\sigma_1 s}\left\|\phi\left(s, a_s\right)-\phi\left(s, b_s\right)\right\|^2 d s\non\\
&\leq e^{-\sigma_1 t} C_\phi^2 \int_{\tau-k}^t e^{\sigma_1 s}\|a(s)-b(s)\|^2 d s \non\\
&\leq e^{-\sigma_1t} C_\phi^2 \int_{\tau-k}^\tau e^{\sigma_1s}\|a(s)-b(s)\|^2 d s+e^{-\sigma_1 t} C_\phi^2 \int_\tau^t e^{\sigma_1s}\|a(s)-b(s)\|^2 d s\non\\
&\leq\frac{C_\phi^2}{\alpha}\left\|\varphi_a-\varphi_b\right\|^2 e^{-\sigma_1(t-\tau)}+\int_\tau^t\|a(s)-b(s)\|^2 d s.
\end{align}

From \eqref{3.58-7}, we obtain
\begin{align}\label{3.59-7}
& \frac{1}{t-\tau}\left(\frac{1}{2 \theta}+\frac{\sigma_2}{2 \lambda_1}\right) e^{-\sigma_1t} \int_\tau^t \int_s^t e^{\sigma_1 r}\left\|\phi\left(r, a_r\right)-\phi\left(r, b_r\right)\right\|^2 d r d s \non\\ & +\frac{1}{t-\tau}\left(\frac{\theta_1}{2 \theta}+\frac{\sigma_3}{2 \lambda_1}\right) e^{-\sigma_1 t} \int_\tau^t e^{\sigma_1s}\left\|\phi\left(s, a_s\right)-\phi\left(s, b_s\right)\right\|^2 d s \non\\ & \leq  \left(\frac{1}{2 \theta}+\frac{\sigma_2}{2 \lambda_1}\right) e^{-\sigma_1t} \int_\tau^t e^{\sigma_1 s}\left\|\phi\left(s, a_s\right)-\phi\left(s, b_s\right)\right\|^2 d s \non\\ & +  \left(\frac{\theta_1}{2 \theta}+\frac{\sigma_3}{2 \lambda_1}\right) \frac{1}{t-\tau} e^{-\sigma_1 t} \int_\tau^t e^{\sigma_1s}\left\|\phi\left(s, a_s\right)-\phi\left(s, b_s\right)\right\|^2 d s\non\\
& \leq \frac{C_\phi^2}{\alpha}\left[\frac{1}{2 \theta}+\frac{\sigma_2}{2 \lambda_1}+\left(\frac{\theta_1}{2 \theta}+\frac{\sigma_3}{2 \lambda_1}\right) \frac{1}{t-\tau}\right]\left\|\varphi_a-\varphi_b\right\|^2 e^{-\sigma_1(t-\tau)} \non\\
& +\frac{C_\phi^2}{\alpha}\left[\frac{1}{2 \theta}+\frac{\sigma_2}{2 \lambda_1}+\left(\frac{\theta_1}{2 \theta}+\frac{\sigma_3}{2 \lambda_1}\right) \frac{1}{t-\tau}\right] \int_\tau^t\|a(s)-b(s)\|^2 d s.
\end{align}

Substituting \eqref{3.59-7} into \eqref{3.57-7} yields
\begin{align}\label{3.60-7}
&I_1(t)+\frac{1-\sigma_1}{t-\tau} e^{-\sigma_1 t} \int_\tau^t e^{\sigma_1 s} I_1(s) d s \non\\
&+\frac{\sigma_2\xi}{t-\tau} e^{-\sigma_1 t} \int_\tau^t \int_s^t e^{\sigma_1 r}\left(\|\nabla a(r)\|^2 \nabla a(r)-\|\nabla b(r)\|^2 \nabla b(r), \nabla v(r)\right) d r d s \non\\
& +\frac{\xi}{t-\tau} e^{-\sigma_1 t} \int_\tau^t \int_s^t e^{\sigma_1 r}\left(\|\nabla a(r)\|^2 \nabla a(r)-\|\nabla b(r)\|^2 \nabla b(r), \nabla \partial_t v(r)\right) d r d s  \non\\
& +\frac{\sigma_3\xi}{t-\tau} e^{-\sigma_1 t} \int_\tau^t e^{\sigma_1 s}\left(\|\nabla a(s)\|^2 \nabla a(s)-\|\nabla b(s)\|^2 \nabla b(s), \nabla v(s)\right) d s  \non\\
& +\frac{\theta_1 \xi}{t-\tau} e^{-\sigma_1 t} \int_\tau^t e^{\sigma_1 s}\left(\|\nabla a(s)\|^2 \nabla a(s)-\| \nabla b(s) \|^2 \nabla b(s), \nabla \partial_t v(s)\right) d s\ \non\\
& +\frac{\sigma_2}{t-\tau} e^{-\sigma_1 t} \int_\tau^t \int_s^t e^{\sigma_1 r}(g(a(r))-g(b(r)), v(r)) d r d s \non\\
& +\frac{1}{t-\tau} e^{-\sigma_1 t} \int_\tau^t \int_s^t e^{\sigma_1 r}(g(a(r))-g(b(r)), \partial_t v(r)) d r d s\non\\
& +\frac{\sigma_3}{t-\tau} e^{-\sigma_1 t} \int_\tau^t e^{\sigma_1 s}(g(a(s))-g(b(s)), v(s)) d s \non\\
& +\frac{\theta_1}{t-\tau} e^{-\sigma_1 t} \int_\tau^t e^{\sigma_1 s}(g(a(s))-g(b(s)), \partial_t v(s)) d s
\non\\
& \leq \frac{\sigma_3}{2(t-\tau)} e^{-\sigma_1(t-\tau)}\|\nabla v(\tau)\|^2+\frac{\sigma_3}{t-\tau} e^{-\sigma_1(t-\tau)}(\partial_t v(\tau), v(\tau))+\frac{\theta_1}{t-\tau} e^{-\sigma_1(t-\tau)} I_1(\tau)\non\\
&-\left(\sigma_2+\frac{\sigma_3}{t-\tau}\right)(\partial_tv(t), v(t))+\frac{\sigma_1 \sigma_3+\sigma_2}{t-\tau} e^{-\sigma_1 t} \int_\tau^t e^{\sigma_1 s}{(\partial_tv(s), v(s)) d s}\non\\
& +\frac{\sigma_1 \sigma_2}{t-\tau} e^{-\sigma_1 t} \int_\tau^t \int_s^t e^{\sigma_1 r}(\partial_t v(r), v(r)) d r d s  \non\\
& +\frac{C_\phi^2}{\alpha}\left[\frac{1}{2 \theta}+\frac{\sigma_2}{2 \lambda_1}+\left(\frac{\theta_1}{2 \theta}+\frac{\sigma_3}{2 \lambda_1}\right) \frac{1}{t-\tau}\right]\left\|\varphi_a-\varphi_b\right\|^2 e^{-\sigma_1(t-\tau)} \non\\
& +\frac{C_\phi^2}{\alpha}\left[\frac{1}{2 \theta}+\frac{\sigma_2}{2 \lambda_1}+\left(\frac{\theta_1}{2 \theta}+\frac{\sigma_3}{2 \lambda_1}\right) \frac{1}{t-\tau}\right] \int_\tau^t\|a(s)-b(s)\|^2 d s \non\\
& +\frac{C_{18}}{t-\tau} e^{-\sigma_1  t}.
\end{align}

Next, we further estimate the terms in \eqref{3.60-7} separately.
Using \eqref{1.2-7} and the embedding $H_0^1(\Omega) \hookrightarrow L^{\frac{2 n}{n-2}}(\Omega)$, we deduce
\begin{align}\label{3.61-7}
& \left|\int_\tau^t e^{\sigma_1 s}(g(a(s))-g(b(s)), v(s)) d s\right| \non\\
& \leq\left(\int_\tau^t \int_\Omega e^{\sigma_1 s}|g(a(s))-g(b(s))|^2 d x d s\right)^{\frac{1}{2}}\left(\int_\tau^t e^{\sigma_1 s}\|v(s)\|^2 d s\right)^{\frac{1}{2}} \non\\
& \leq C_{19} e^{\frac{\sigma_1}{2} t}\left(\frac{|\Omega|}{\sigma_1} e^{\sigma_1 t}+\left(\int_\tau^t e^{\sigma_1 s}\left(\|\nabla a(s)\|^2+\|\nabla b(s)\|^2\right) d s\right)^{\frac{n}{n-2}}\right)^{\frac{1}{2}}\left(\int_\tau^t\|\nabla v(s)\|^2 d s\right)^{\frac{1}{2}}.
\end{align}

From Lemma \ref{lem3.2-7} and \eqref{3.61-7}, we obtain there exists a constant $C_{20}=C_{20}\left(B_0\right)>0$ such that
\begin{align}\label{3.62-7}
& \frac{\sigma_3}{t-\tau} e^{-\sigma_1 t} \int_\tau^t e^{\sigma_1s}(g(a(s))-g(b(s)), v(s)) d s \non\\
& \left.\geq-\frac{\sigma_3 C_{19}}{t-\tau} e^{-\frac{\sigma_1}{2} t}\left(\frac{|\Omega |}{\sigma_1} e^{\sigma_1 t}+\int_\tau^t e^{\sigma_1 s}\left(\|\nabla a(s)\|^2+\|\nabla b(s)\|^2\right) d s\right)^{\frac{n}{n-2}}\right)^{\frac{1}{2}}\left(\int_\tau^t\|\nabla v(s)\|^2 d s\right)^{\frac{1}{2}}
\non\\
&\geq-\frac{\sigma_3 C_{20}}{t-\tau}\left(\int_\tau^t\|\nabla v(s)\|^2 d s\right)^{\frac{1}{2}}.
\end{align}

Moreover, there exists a constant $C_{21}=C_{21}(B_{0})>0$ such that
\begin{align}\label{3.63-7}
& \frac{\theta_1}{t-\tau} e^{-\sigma_1 t} \int_\tau^t e^{\sigma_1 s}(g(a(s))-g(b(s)), \partial_t v(s) d s \geq - \frac{\theta_1C_{21}}{t-\tau}\left(\int_\tau^t \| \partial_t v(s)\|^2 d s\right)^{\frac{1}{2}} .
\end{align}

Using the Sobolev embedding $H_0^1(\Omega) \hookrightarrow L^{\frac{2n}{n-2}}(\Omega)$,
together with \eqref{1.2-7}, the H\"older inequality, and Lemma \ref{lem3.2-7},
we infer that there exists a constant $C_{22}=C_{22}(B_0)>0$ such that
\begin{align}\label{3.64-7}
&\left|\frac{\sigma_2}{t-\tau}\mathrm{e}^{-\sigma_1 t} \int_\tau^t \int_s^t e^{\sigma_1 r}(g(a(r))-g(b(r)), v(r)) d r d s\right|\non\\
& \leq \sigma_2 e^{-\sigma_1 t}\left|\int_s^t e^{\sigma_1 r}(g(a(r))-g(b(r)), v(r))d r\right|\non \\
& \leq\sigma_2 e^{-\sigma_1 t}\left|\int_\tau^t e^{\sigma_1 s}(g(a(s))-g(b(s)), v(s))d s\right|\non \\
& \leq \sigma_2 e^{-\sigma_1 t}\left(\int_\tau^t e^{\sigma_1 s} \int_{\Omega}|g(a(s))-g(b(s))|^2 d x d s\right)^{\frac{1}{2}}\left(\int_\tau^t e^{\sigma_1 s}\|v(s)\|^2 d s\right)^{\frac{1}{2}}\non\\
& \leq \sigma_2 e^{-\frac{\sigma_1}{2} t}\left(\frac{|\Omega|}{\sigma_1} e^{\sigma_1 t}+\left(\int_\tau^t e^{\sigma_1 s}\left(\|\nabla a(s)\|^2+\|\nabla b(s)\|^2\right) d s\right)^{\frac{n}{n-2}}\right)^{\frac{1}{2}}\left(\int_\tau^t\|v(s)\|^2 d s\right)^{\frac{1}{2}}\non\\
& \leq \sigma_2 C_{22}\left(\int_\tau^t\|v(s)\|^2 d s\right)^{\frac{1}{2}}.
\end{align}

Then we obtain there exists a constant $C_{23}>0$ such that
\begin{align}\label{3.65-7}
&\frac{\sigma_2}{t-\tau} e^{-\sigma_1 t} \int_\tau^t \int_s^t e^{\sigma_1 r}(g(a(r))-g(b(r)), v(r)) d r d s \geq-\sigma_2 C_{23}\left(\int_\tau^t\|v(s)\|^2 d s\right)^{\frac{1}{2}}.
\end{align}

Similar to \eqref{3.65-7}, we deduce there exists a constant $C_{24}>0$ such that
\begin{align}\label{3.66-7}
& \frac{1}{t-\tau} e^{-\sigma_1t} \int_\tau^t \int_s^t e^{\sigma_1 r}(g(a(r))-g(b(r)), \partial_t v(r)) d r d s \geq-C_{24}\left(\int_\tau^t\|\partial_t v(s)\|^2 d s\right)^{\frac{1}{2}}.
\end{align}

Using Lemma \ref{lem3.2-7} and the H\"older inequality, we deduce there exists a constant $C_{25}=C_{25}\left(B_0\right)>0$ such that
\begin{align}\label{3.67-7}
& \frac{\sigma_1 \sigma_3+\sigma_2}{t-\tau} e^{-\sigma_1 t} \int_\tau^t e^{\sigma_1 s}(\partial_t v(s), v(s)) d s \non\\
& \leq \frac{\sigma_1 \sigma_3+\sigma_2}{t-\tau}\left(e^{-\sigma_1 t} \int_\tau^t e^{\sigma_1 s}\|\partial_t v(s)\|^2 d s\right)^{\frac{1}{2}}\left(e^{-\sigma_1 t} \int_\tau^t e^{\sigma_1 s}\|v(s)\|^2 d s\right)^{\frac{1}{2}} \non\\
& \leq \frac{\left(\sigma_1 \sigma_3+\sigma_2\right) C_{25}}{t-\tau}\left(\int_\tau^t\|v(s)\|^2 d s\right)^{\frac{1}{2}}.
\end{align}

Similar to \eqref{3.67-7}, we derive there exists a constant $C_{26}=C_{26}\left(B_0\right)>0$ such that
\begin{align}\label{3.68-7}
&\frac{\sigma_1 \sigma_2}{t-\tau} e^{-\sigma_1 t}\left|\int_\tau^t \int_s^t e^{\sigma_1 r}\left(\partial_t v(r), v(r)\right) d r d s\right|\non\\
& \leq\sigma_1 \sigma_2 e^{-\sigma_1 t}\left|\int_s^t e^{\sigma_1 r}\left(\partial_t v(r), v(r)\right) d r\right|\non \\
&\leq \sigma_1 \sigma_2 e^{-\sigma_1 t}\left(e^{-\sigma_1 t} \int_\tau^t e^{\sigma_1 s}\|\partial_t v(s)\|^2 d s\right)^{\frac{1}{2}}\left(e^{-\sigma_1 t} \int_\tau^t e^{\sigma_1 s}\|v(s)\|^2 d s\right)^{\frac{1}{2}}\non\\
&\leq \sigma_1 \sigma_2 C_{26}\left(\int_\tau^t\|v(s)\|^2 d
s\right)^{\frac{1}{2}}.
\end{align}

By the Poincar\'e inequality and noting that $v=a-b$, we conclude
\begin{align}\label{3.69-7}
& \frac{\sigma_2\xi}{t-\tau} e^{-\sigma_1 t} \int_\tau^t \int_s^t e^{\sigma_1 r}\left(\|\nabla a(r)\|^2 \nabla a(r)-\|\nabla b(r)\|^2 \nabla b(r), \nabla v(r)\right) d r d s \non\\
& +\frac{\sigma_3\xi}{t-\tau} e^{-\sigma_1 t} \int_\tau^t e^{\sigma_1 s}\left(\|\nabla a(s)\|^2 \nabla a(s)-\|\nabla b(s)\|^2 \nabla b(s), \nabla v(s)) d s\right. \non\\
& =\frac{\sigma_2\xi}{t-\tau} e^{-\sigma_1 t} \int_\tau^t \int_s^t e^{\sigma_1 r}\left(\|\nabla a(r)\|^2 \nabla v(r), \nabla v(r)\right) d r d s \non\\
& +\frac{\sigma_2\xi}{t-\tau} e^{-\sigma_1t} \int_\tau^t \int_s^t e^{\sigma_1r}\left(\left(\|\nabla a(r)\|^2-\|\nabla b(r)\|^2\right) \nabla b(r), \nabla v(r)\right) d r d s \non\\
& +\frac{\sigma_3\xi}{t-\tau} e^{-\sigma_1 t} \int_\tau^t e^{\sigma_1 s}\left(\|\nabla a(s)\|^2 \nabla v(s), \nabla v(s)\right) d s\non\\
& +\frac{\sigma_3\xi}{t-\tau} e^{-\sigma_1 t }\int_\tau^t e^{\sigma_1 s}\left(\left(\|\nabla a(s) \|^2- \|\nabla b(s) \|^2\right) \nabla b(s), \nabla v(s)\right) d s\non\\
& \equiv V_1(t)
\end{align}
and
\begin{align}\label{3.70-7}
& \frac{\xi\lambda_1}{t-\tau} e^{-\sigma_1 t} \int_\tau^t \int_s^t e^{\sigma_1 r}\left(\|\nabla a(r)\|^2 \nabla a(r)-\|\nabla b(r)\|^2 \nabla b(r), \partial_t v(r)\right) d r d s \non\\
& +\frac{\theta_1 \xi \lambda_1}{t-\tau}e^{-\sigma_1 t} \int_\tau^t e^{\sigma_1 s}\left(\|\nabla a(s)\|^2 \nabla a(s)-\|\nabla b(s)\|^2 \nabla b(s), \partial_t v(s)\right) d s \non\\
& =\frac{\xi \lambda_1}{t-\tau} e^{-\sigma_1 t} \int_\tau^t \int_s^t e^{\sigma_1 r}\left(\|\nabla a(r)\|^2 \nabla v(r), \partial_t v(r)\right) d r d s \non\\
& +\frac{\xi \lambda_1}{t-\tau} e^{-\sigma_1 t} \int_\tau^t \int_s^t e^{\sigma_1 r}\left(\|\nabla a(r)\|^2-\|\nabla b(r)\|^2\right) \nabla b(r), \partial_t v(r)) d r d s \non\\
& +\frac{\theta_1 \xi \lambda_1}{t-\tau}e^{-\sigma_1 t} \int_\tau^t e^{\sigma_1 s}\left(\|\nabla a(s)\|^2 \nabla v(s), \partial_tv(s)\right) d s \non\\
& +\frac{\theta_1 \xi \lambda_1}{t-\tau} e^{-\sigma_1 t} \int_\tau^t e^{\sigma_1 s}\left(\left(\|\nabla a(s)\|^2-\|\nabla b(s)\|^2\right) \nabla b(s), \partial_t v(s)\right) d s \non\\
& \equiv V_2(t).
\end{align}

$\text{Substituting \eqref{3.62-7}, \eqref{3.63-7}, and \eqref{3.65-7}$-$\eqref{3.70-7} into \eqref{3.60-7}, and noting that } \sigma_1<1,\text{ we obtain}$
\begin{align}\label{3.71-7}
&I_1(t) \leq \frac{\sigma_3}{2(t-\tau)} e^{-\sigma_1(t-\tau)}\|\nabla v(\tau)\|^2+\frac{\sigma_3}{t-\tau} e^{-\sigma_1(t-\tau)}(\partial_tv(\tau), v(\tau))\non\\
&+\frac{\theta_1}{t-\tau} e^{-\sigma_1(t-\tau)} I_1(\tau) -\left(\sigma_2+\frac{\sigma_3}{t-\tau}\right)(\partial_tv(t), v(t))\non\\
& +\frac{C_{18}}{t-\tau} e^{-\sigma_1 t}-V_1(t)-V_2(t)+\frac{\sigma_3 C_{20}}{t-\tau}\left(\int_\tau^t\|\nabla v(s)\|^2 d s\right)^{\frac{1}{2}} \non\\
&+C_{27}\left(\int_\tau^t\|\partial_t v(s)\|^2 d s\right)^{\frac{1}{2}}+C_{28}\left(\int_\tau^t\|v(s)\|^2 d s\right)^{\frac{1}{2}} \non\\
& +\frac{C_\phi^2}{\alpha}\left[\frac{1}{2 \theta}+\frac{\sigma_2}{2 \lambda_1}+\left(\frac{\theta_1}{2 \theta}+\frac{\sigma_3}{2 \lambda_1}\right) \frac{1}{t-\tau}\right]\left\|\varphi_a-\varphi_b\right\|^2 e^{-\sigma_1(t-\tau)} \non\\
& +\frac{C_\phi^2}{\alpha}\left[\frac{1}{2 \theta}+\frac{\sigma_2}{2 \lambda_1}+\left(\frac{\theta_1}{2 \theta}+\frac{\sigma_3}{2 \lambda_1}\right) \frac{1}{t-\tau}\right] \int_\tau^t\|a(s)-b(s)\|^2 d s.
\end{align}
where $
C_{27}=\frac{\theta_1 C_{21}}{t-\tau}+C_{24}$ and $
C_{28}=\sigma_2 C_{23}+\frac{\left(\sigma_1 \sigma_3+\sigma_2\right) C_{25}}{t-\tau}+\sigma_1 \sigma_2 C_{26}$.

Furthermore, let
\begin{align}\label{3.72-7}
Y(t)&=\frac{\sigma_3}{2(t-\tau)} e^{-\sigma_1(t-\tau)}\|\nabla v(\tau)\|^2+\frac{\sigma_3}{t-\tau} e^{-\sigma_1(t-\tau)}(\partial_tv(\tau), v(\tau))+\frac{\theta_1}{t-\tau} e^{-\sigma_1(t-\tau)} I_1(\tau) \non\\
& +\frac{C_\phi^2}{\alpha}\left[\frac{1}{2 \theta}+\frac{\sigma_2}{2 \lambda_1}+\left(\frac{\theta_1}{2 \theta}+\frac{\sigma_3}{2 \lambda_1}\right) \frac{1}{t-\tau}\right]\left\|\varphi_a-\varphi_b\right\|^2 e^{-\sigma_1(t-\tau)}+\frac{C_{18}}{t-\tau} e^{-\sigma_1 t}
\end{align}
and
\begin{align}\label{3.73-7}
&\psi(\eta_a,\eta_b)=\psi((a,\partial_t a),(b,\partial_t b))\non\\
&=-\left(\sigma_2+\frac{\sigma_3}{t-\tau}\right)(\partial_t a(t)-\partial_t b(t), a(t)-b(t))+\frac{\sigma_3 C_{20}}{t-\tau}\left(\int_\tau^t\|\nabla(a(s)-b(s))\|^2 d s\right)^{\frac{1}{2}} \non\\
& +C_{27}\left(\int_\tau^t(\partial_t a(s)-\partial_t b(s)) d s\right)^{\frac{1}{2}} +C_{28}\left(\int_\tau^t\|a(s)-b(s)\|^2 d s\right)^{\frac{1}{2}} \non\\
& +\frac{C_{\phi}^2}{\alpha}\left[\frac{1}{2 \theta}+\frac{\sigma_2}{2 \lambda_1}+\left(\frac{\theta_1}{2 \theta}+\frac{\sigma_3}{2 \lambda_1}\right) \frac{1}{t-\tau}\right] \int_\tau^t\|a(s)-b(s)\|^2 d s\non \\
& -V_1(t)-V_2(t).
\end{align}

Then from \eqref{3.71-7}$-$\eqref{3.73-7}, we derive there exists a constant $\varepsilon>0$ such that
\begin{align}\label{3.74-7}
I_1(t) \leq \varepsilon+\psi((a, \partial_ta),(b, \partial_t b)).
\end{align}

Therefore, by setting $T=t-\tau$ and applying Lemma \ref{Lem1-7}, we conclude that the asymptotic compactness of $\{U(\cdot,\cdot)\}_{t\geq \tau}$  follows provided that
\begin{equation}\label{3.75-7}
\psi((a,\partial_t a),(b,\partial_t b)) \in C(B_0).
\end{equation}

Let $\left(v^i,\partial_t v^i\right)$ be the solution of equation \eqref{1.1-7} corresponding to the initial data $\left(\varphi^i,\partial_t\varphi^i\right)\in B_0(\tau)\times B_0(\tau)$, where $i\in\mathbb{N}^+$. Then, by Theorem \ref{lem3.0-7}, it follows that
\begin{equation}\label{3.76-7}
v^i \rightarrow v\quad \text{weakly-star in }
L^{\infty}\left(\tau, t ; H_0^1(\Omega)\right),
\end{equation}
\begin{equation}\label{3.77-7}
\partial_t v^i \rightarrow \partial_t v
\quad \text{weakly-star in }
L^{\infty}\left(\tau, t ; L^2(\Omega)\right),
\end{equation}
\begin{equation}\label{3.78-7}
v^i \to v\quad \text{in }
L^2(\Omega) \text{ and } L^{\frac{2n-2}{n-2}}(\Omega),
\end{equation}
\begin{equation}\label{3.79-7}
v^i \to v\quad \text{in }
L^2\left(\tau, t ; L^2(\Omega)\right).
\end{equation}

By \eqref{3.74-7}$-$\eqref{3.77-7} and by the Aubin-Lions lemma, it follows that there exists a subsequence $\left(v^j, \partial_tv^j\right)$ of $\left(v^i, \partial_tv^i\right)$ such that
\begin{equation}
\lim \limits_{i \rightarrow+\infty} \lim\limits _{j \rightarrow+\infty} \int_{\Omega}\left(\partial_t v^i-\partial_t v^j\right)\left(v^i-v^j\right) d x=0,
\label{3.80-7}
\end{equation}
\begin{equation}
\lim\limits_{i \rightarrow+\infty} \lim \limits_{j \rightarrow+\infty} \int_\tau^t\left\|\nabla v^i-\nabla v^j\right\|^2 d s=0,
\label{3.81-7}
\end{equation}
\begin{equation}
\lim\limits _{i \rightarrow+\infty} \lim\limits _{j \rightarrow+\infty} \int_\tau^t\left(\partial_t v^i-\partial_t v^j\right) d s=0,
\label{3.82-7}
\end{equation}
\begin{equation}
\lim\limits _{i \rightarrow+\infty} \lim\limits _{j \rightarrow+\infty} \int_\tau^t\left\|v^i-v^j\right\|^2 d s=0,
\label{3.83-7}
\end{equation}
\begin{align}\label{3.84-7}
& \lim _{i \rightarrow+\infty} \lim _{j \rightarrow+\infty} \int_\tau^t \int_s^t e^{\sigma_1 r}\left(\left\|\nabla v^i(r)\right\|^2 \nabla v^i(r)-\left\|\nabla v^j(r)\right\|^2 \nabla v^j(r), \nabla\left(v^i(r)-v^j(r)\right)\right) d r d s \non\\
& =\lim _{i \rightarrow+\infty} \lim _{j \rightarrow+\infty} \int_\tau^t \int_s^t e^{\sigma_1 r}\left(\left\|\nabla v^i(r)\right\|^2 \nabla\left(v^i(r)-v^j(r)\right), \nabla\left(v^i(r)-v^j(r)\right)\right) d r d s \non\\
& +\lim _{i \rightarrow+\infty} \lim _{j \rightarrow+\infty} \int_\tau^t \int_s^t e^{\sigma_1 r}\left(\left(\left\|\nabla v^i(r)\right\|^2-\left\|\nabla v^j(r)\right\|^2\right) \nabla v^j(r), \nabla\left(v^i(r)-v^j(r)\right)\right) d r d s \non\\
& =0,
\end{align}
\begin{align}\label{3.85-7}
& \lim _{i \rightarrow+\infty} \lim _{j \rightarrow+\infty} \int_\tau^t e^{\sigma_1 s}\left(\left\|\nabla v^i(s)\right\|^2 \nabla v^i(s)-\left\|\nabla v^j(s)\right\|^2 \nabla v^j(s), \nabla\left(v^i(r)-v^j(r)\right)\right) d s \non\\
& =\lim _{i \rightarrow+\infty} \lim _{j \rightarrow+\infty} \int_\tau^t e^{\sigma_1 s}\left(\left\|\nabla v^i(s)\right\|^2 \nabla\left(v^i(s)-v^j(s)\right), \nabla\left(v^i(s)-v^j(s)\right)\right) d s \non\\
& +\lim _{i \rightarrow+\infty} \lim _{j \rightarrow+\infty} \int_\tau^t e^{\sigma_1 s}\left(\left(\left\|\nabla v^i(s)\right\|^2-\left\|\nabla v^j(s)\right\|^2\right) \nabla v^j(s), \nabla\left(v^i(s)-v^j(s)\right)\right) d s \non\\
& =0,
\end{align}
\begin{align}\label{3.86-7}
& \lim _{i \rightarrow+\infty} \lim _{j \rightarrow+\infty} \int_\tau^t \int_s^t e^{\sigma_1 r}\left(\left\|\nabla v^i(r)\right\|^2 \nabla v^i(r)-\left\|\nabla v^j(r)\right\|^2 \nabla v^j(r), \partial_t v^i(r)-\partial_t v^j(r)\right) d r d s \non\\
& =\lim _{i \rightarrow+\infty} \lim _{j \rightarrow+\infty} \int_\tau^t \int_s^t e^{\sigma_1 r}\left(\left\|\nabla v^i(r)\right\|^2 \nabla\left(v^i(r)-v^j(r)\right), \partial_t v^i(r)-\partial_t v^j(r)\right) d r d s \non\\
&+\lim _{i \rightarrow+\infty} \lim _{j \rightarrow+\infty} \int_\tau^t \int_s^t e^{\sigma_1 r}\left(\left(\left\|\nabla v^i(r)\right\|^2-\left\|\nabla v^j(r)\right\|^2\right) \nabla v^j(r), \partial _t v^i(v)-\partial_t v^j(r)\right) d r d s\non\\
& =0,
\end{align}
\begin{align}\label{3.87-7}
& \lim _{i \rightarrow+\infty} \lim _{j \rightarrow+\infty} \int_\tau^t e^{\sigma_1 s}\left(\left\|\nabla v^i(s)\right\|^2 \nabla v^i(s)-\left\|\nabla v^j(s)\right\|^2 \nabla v^j(s), \partial_t v^i(s)-\partial_t v^j(s)\right) d s \non\\
& =\lim _{i \rightarrow+\infty} \lim _{j \rightarrow+\infty} \int_\tau^t e^{\sigma_1 s}\left(\left\|\nabla v^i(s)\right\|^2 \nabla\left(v^i(s)-v^j(s)\right), \partial_t v^i(s)-\partial_t v^j(s)\right) d s \non\\
& +\lim _{i \rightarrow+\infty} \lim _{j \rightarrow+\infty} \int_\tau^t e^{\sigma_1 s}\left(\left(\left\|\nabla v^i(s)\right\|^2-\left\|\nabla v^j(s)\right\|^2\right) \nabla v^j(s), \partial_t v^i(s)-\partial_t v^j(s)\right) d s \non\\
& =0.
\end{align}

It follows from \eqref{3.71-7}$-$\eqref{3.73-7} and \eqref{3.80-7}$-$\eqref{3.87-7} that
\(\psi\left(\eta_a,\eta_b\right)\) is a contractive function on
\(B_0(\tau)\times B_0(\tau)\). Consequently, the process
\(\{v(t,\tau)\}_{t\geq \tau}\) is pullback asymptotically compact in
\(C_{H_0^{1}(\Omega)}\times C_{L^2(\Omega)}\).

$\hfill$$\Box$

\begin{theorem}\label{thm4.6}
Under the assumptions of Lemma \rm{\ref{lem3.3-7}}, the process
$\{U(t,\tau)\}_{t\geq \tau}$ possesses a pullback
$\mathcal{D}$-attractor in
$C_{H_0^1(\Omega)} \times C_{L^2(\Omega)}$.
\end{theorem}
\noindent $\mathbf{Proof}$ From Lemmas \ref{LemmaLH}, \ref{lem3.2-7}, and \ref{lem3.3-7}, Theorem \ref{thm4.6} follows directly.

$\hfill$$\Box$

\noindent$\mathbf{Acknowledgment}$

We would like to express our sincere gratitude to the editors and anonymous reviewers for their invaluable feedback and insightful comments on the initial draft of this paper. Their thorough review and constructive suggestions have played a crucial role in improving the quality of this manuscript. We greatly appreciate their dedication and support throughout the review process.

\noindent$\mathbf{Funding}$

This work was supported by the National Natural Science Foundation of China (Grant No. 12462005; 62063025), the Natural Science Foundation of Inner Mongolia Autonomous Region of China (Grant No. 2024MS06028; 2026QC0419; 2026MS0464; 2026MS0469; 2026MS0478) and the Keju Plan of Inner Mongolia University of Science and Technology (Grant No. KJJH2024985; KJJH2025991).

\noindent$\mathbf{Conflict\,\,of\,\,interest\,\,statement}$

The authors have no conflict of interest.

\end{document}